\documentclass[runningheads,12pt]{llncs}

\usepackage{amssymb}
\usepackage [english]{babel}
\usepackage[utf8]{inputenc}
\usepackage[T1]{fontenc}

\usepackage{amssymb}
\usepackage{amsmath}
\usepackage{amsfonts}

\usepackage{graphicx}
\usepackage{url}

\usepackage[usenames]{color}

\usepackage{algorithm}
\usepackage{algpseudocode}

\usepackage{parskip}
\usepackage{subcaption}
\usepackage{makecell}

\usepackage{xifthen}

  \usepackage{geometry}
 \usepackage{lineno}

\newcommand{\PathToPic}{}

\DeclareMathOperator{\mdiv}{div}
\DeclareMathOperator{\grad}{\nabla}

\newcommand{\p}[1]{\hat{#1}}
\newcommand{\g}[1]{#1}

\newcommand{\sMP}[1][]{\ifthenelse{\isempty{#1}}{^{(k)}}{^{(#1)}}}

\newcommand{\ifSet}{\iSet_{\mathcal{F}}}
\newcommand{\iSet}{{\mathcal{I}}}

 \newcommand{\MatTwo}[4]
  {\begin{bmatrix}
    #1 & #2\\
    #3 & #4
  \end{bmatrix}}
  
  \newcommand{\VecTwo}[2]
  {\begin{bmatrix}
    #1 \\
    #2 
  \end{bmatrix}}
  
  \newcommand{\MatThree}[9]
  {
    \begin{bmatrix}
    #1 & #2 & #3\\
    #4 & #5 & #6\\
    #7 & #8 & #9
    \end{bmatrix}
  }
  
  \newcommand{\VecThree}[3]
  {\begin{bmatrix}
    #1 \\
    #2 \\
    #3
  \end{bmatrix}}

    \newcommand{\keywords}[1]{\par\addvspace\baselineskip  
     \noindent\keywordname\enspace\ignorespaces#1}

\begin{document}

\urldef{\mailsa}\path|christoph.hofer@jku.at|

\title{Parallelization of continuous and discontinuous Galerkin dual-primal Isogeometric tearing and interconnecting methods}
\author{Christoph Hofer$^1$}
\institute{ $^1$ Johannes Kepler University (JKU),
				Altenbergerstr. 69, A-4040 Linz, Austria,\\
\mailsa 
 }

\noindent
\maketitle







\begin{abstract}
In this paper we investigate the parallelization of dual-primal isogeometric tearing and interconnecting (IETI-DP) type methods for solving large-scale continuous and discontinuous Galerkin systems of equations arising from Isogeometric analysis of elliptic boundary value problems.
These methods are extensions of the finite element tearing and interconnecting methods to isogeometric analysis. The algorithms are implemented by means of energy minimizing primal subspaces. We discuss how these methods can efficiently be parallelized in a distributed memory setting. Weak and strong scaling studies presented for two and three dimensional problems show an excellent parallel efficiency.
\end{abstract}



\pagestyle{myheadings}
\thispagestyle{plain}
\markboth{}{C. Hofer, Parallelization of cG and dG-IETI-DP methods}

\section{Introduction}
Isogeometric Analysis (IgA) is a novel methodology for the numerical solution of partial differential equations (PDE). IgA was first introduced by Hughes, Cottrell and Bazilevs in \cite{HL:HughesCottrellBazilevs:2005a}, see also the monograph \cite{HL:CotrellHughesBazilevs:2009a} for a comprehensive presentation of the IgA framework and the recent survey article \cite{HL:BeiraodaVeigaBuffaSangalliVazquez:2014a}. The main principle is to use the same basis functions for describing the geometry and to represent the discrete solution of the PDE problem under consideration. The most common choices are B-Splines, Non Uniform Rational B-Splines (NURBS), T-Splines, Truncated Hierarchical B-Splines (THB-Splines), etc., see, e.g., \cite{HL:GiannelliJuettlerSpeleers:2012a}, \cite{HL:GiannelliJuettlerSpeleers:2014a} and  \cite{HL:BazilevsCaloCottrellEvans:2010a}.
One of the strengths of IgA is the capability of creating high-order splines spaces, while keeping the number of degrees of freedom quite small. 
Moreover, having basis functions with high smoothness is useful when considering higher-order PDEs, e.g., the biharmonic equation.

In many cases the domain can not be represented with a single mapping, referred to as \emph{geometrical mapping}.  Complicated geometries are decomposed into simple domains, called \emph{patches}, which are topologically equivalent to a cube. The set of patches forming the computational domain is called multipatch domain.
The obtained patch parametrizations and the original geometry may not be identical. The result are
small gaps and overlaps  occurring  at the interfaces of the patches, called \emph{segmentation crimes}, see \cite{HL:JuettlerKaplNguyenPanPauley:2014a},
 \cite{HL:PauleyNguyenMayerSpehWeegerJuettler:2015a} and 
 \cite{Hoschek_Lasser_CAD_book_1993} for a comprehensive analysis. Nevertheless, one still wants to solve PDEs on such domains. To do so, numerical schemes based on the discontinuous Galerkin (dG) method for elliptic PDEs
 were developed in \cite{HL:HoferLangerToulopoulos:2016a}, \cite{HL:HoferToulopoulos:2016a} and \cite{HL:HoferLangerToulopoulos:2016b}. There, the corresponding error analysis is also provided. In addition to domains with segmentation crimes, the dG formulation is very useful when considering different B-Splines spaces on each patch, e.g., non-matching grids at the interface and different spline degrees. An analysis for the dG-IgA formulation with extensions to low regularity solutions can be found in \cite{HL:LangerToulopoulos:2015a}. For a detailed discussion of dG for finite element methods, we refer, e.g., to \cite{HL:Riviere:2008a} and \cite{HL:PietroErn:2012a}. 
 
 In the present paper, we are considering fast solution methods for linear systems arising from the discretization of elliptic PDEs by means of IgA.
 We investigate non-overlapping domain decomposition (DD) methods of the dual-primal tearing and interconnecting type. This type of methods are equivalent to the so called Balancing Domain Decomposition by Constraints (BDDC) methods, see \cite{HL:MandelDohrmannTezaur:2005a}, \cite{HL:ToselliWidlund:2005a}, \cite{HL:Pechstein:2013a} and references therein. The version based on a conforming Galerkin (cG) discretization, called  dual-primal isogeometric tearing and interconnecting (IETI-DP) method was first introduced in \cite{HL:KleissPechsteinJuettlerTomar:2012a} and the equivalent IgA BDDC method was analysed in \cite{HL:VeigaChoPavarinoScacchi:2013a}. Further extensions to the analysis are presented in \cite{HL:HoferLanger:2016b}. The version based on the dG formulation, abbreviated by dG-IETI-DP, was introduced in \cite{HL:HoferLanger:2016a} and analyzed in \cite{HL:Hofer:2016a}, see \cite{HL:DryjaGalvisSarkis:2007a}, \cite{HL:DryjaGalvisSarkis:2013a} and \cite{HL:DryjaSarkis:2014a} for the corresponding finite element counterparts. We also want to mention development in overlapping Schwarz methods, see, e.g., \cite{HL:VeigaChoPavarinoScacchi:2012a} and \cite{HL:VeigaChoPavarinoScacchi:2013b}.  The aim of this paper is to present the parallel scalability of the  cG and dG IETI-DP methods. We investigate weak and strong scaling in two and three dimensional domains for different B-Spline degrees. The implemented algorithms are based on energy minimizing primal subspaces, which simplifies the parallelization of the solver part, but having more effort in the setup phase (assembling phase). We rephrase key parts of this algorithm and discuss how to realize the communication by means of Message Passing Interface (MPI). 
 In general, FETI-DP and equivalent BDDC methods are by nature well suited for large-scale parallelization and has been widely studied for solving large-scale finite element equations, e.g.,  in \cite{HL:KlawonnRheinbach:2010a}, \cite{HL:Rheinbach:2009}, \cite{HL:KlawonnRheinbach:2006a}
 and \cite{HL:KlawonnLanserRheinbach:2015a}, see also \cite{HL:KlawonnLanserRheinbachStengelWellein:2015a} for a hybrid OpenMP/MPI version. Considering  a domain decomposition with several ten thousands of subdomains, the influence of the coarse grid problem becomes more and more significant. Especially, its LU-factorization is the bottleneck of the algorithm. The remedy is to reformulate the FETI-DP system in such a way that the solution of the coarse grid problem is not required in the application of the system matrix, but in the preconditioner. This enables the use of inexact methods like geometric or algebraic multigrid, see, e.g., \cite{HL:KlawonnLanserRheinbach:2016a}, \cite{HL:KlawonnLanserRheinbach:2015a}, \cite{HL:KlawonnRheinbach:2007a},  \cite{HL:KlawonnRheinbach:2010a} and \cite{HL:KlawonnRheinbachPavarino:2008a}. Moreover, inexact solvers can also be used in the scaled Dirichlet preconditioner and, if using the saddle point formulation, also for the local solvers, cf., \cite{HL:KlawonnRheinbach:2007a}, see also \cite{HL:KlawonnRheinbach:2010a}, \cite{HL:Rheinbach:2009} and references therein for alternative approaches by means of hybrid FETI. We also want to mention inexact version for the BDDC method, see, e.g., \cite{HL:Tu:2007b}, \cite{HL:Tu:2007a}, \cite{HL:Dohrmann:2007a}, \cite{HL:LiWidlund:2007a} and \cite{HL:Zampini:2014}. FETI-DP methods has also been successfully applied to non-linear problems my means of a non-linear version of FETI-DP. We want to highlight recent advances presented, e.g., in \cite{HL:KlawonnLanserRheinbach:2014a}, \cite{HL:KlawonnLanserRheinbach:2016a} and \cite{HL:KlawonnLanserRheinbach:2015a}, showing excellent scalability on large-scale supercomputers. 
    
 In the present paper, we consider the following 
second-order elliptic boundary value problem 
in a bounded Lipschitz domain $\Omega\subset \mathbb{R}^d,$ with $d\in\{2,3\}$:
Find $u: \overline{\Omega} \rightarrow \mathbb{R}$ such that\\ 
\begin{equation}
  \label{equ:ModelStrong}
  - \mdiv(\alpha \grad u) = f \; \text{in } \Omega,\;
  u = 0 \; \text{on } \Gamma_D,  \;\text{and}\;
  \alpha \frac{\partial u}{\partial n} = g_N \; \text{on } \Gamma_N,
\end{equation}
with given, sufficient smooth data $f, g_N \text{ and } \alpha$, where  
the coefficient $\alpha$ is uniformly bounded from below and above
by some positive constants $\alpha_{min}$ and $\alpha_{max}$, respectively.
The boundary 
$\partial \Omega$ 
of the computational domain $\Omega$
consists of a Dirichlet part $\Gamma_D$ of positive boundary measure 
and a Neumann part $\Gamma_N$.
Furthermore, we assume that the Dirichlet boundary $\Gamma_D$ is always 
a union of complete patch sides (edges / face in 2D / 3D) 
which are uniquely defined in IgA.
Without loss of generality, we assume 
homogeneous
Dirichlet conditions. 
This can always be  obtained by homogenization.
By means of integration by parts, 
we arrive at
the weak formulation of \eqref{equ:ModelStrong} 
which reads as follows:
Find $u \in V_{D} = \{ u\in H^1: \gamma_0 u = 0 \text{ on } \Gamma_D \}$
such that
\begin{align}
  \label{equ:ModelVar}
    a(u,v) = \left\langle F, v \right\rangle \quad \forall v \in V_{D},
\end{align}
where $\gamma_0$ denotes the trace operator. The bilinear form
$a(\cdot,\cdot): V_{D} \times V_{D} \rightarrow \mathbb{R}$
and the linear form $\left\langle F, \cdot \right\rangle: V_{D} \rightarrow \mathbb{R}$
are given by the expressions
\begin{equation*}
a(u,v) := \int_\Omega \alpha \nabla u \cdot \nabla v \,dx
\quad \mbox{and} \quad
\left\langle F, v \right\rangle := \int_\Omega f v \,dx + \int_{\Gamma_N} g_N v \,ds.
\end{equation*}
The remainder of the paper is organized as follows: In Section~\ref{sec:iga}, we give a short introduction to isogeometric analysis, providing the basic definitions and notations. Section~\ref{sec:galerkin-IGA} describes the different discretizations of the model problem obtained the continuous and discontinuous Galerkin methods. In Section~\ref{sec:IETI-method}, we formulate the IETI-DP method for both discretizations and provide implementation details. The way how the algorithm is parallelized is explained in Section~\ref{sec:para}. Numerical examples are presented in Section~\ref{sec:num_ex}. Finally we draw some conclusions in Section~\ref{sec:conclusion}.

\section{Isogeometric Analysis}
\label{sec:iga}
In this section, we give a very short overview about IgA. For a more comprehensive study, we refer to, e.g., \cite{HL:CotrellHughesBazilevs:2009a} and \cite{HL:LangerToulopoulos:2015a}.

Let $\p{\Omega}:=(0,1)^d,d\in\{2,3\}$, be the d-dimensional unit cube, which we refer to as the \emph{parameter domain}. Let $p_\iota$ and $M_\iota,\iota\in\{1,\ldots,d\}$, be the B-Spline degree and the number of basis functions in $x_\iota$-direction. Moreover, let $\Xi_\iota = \{\xi_1=0,\xi_2,\ldots,\xi_{n_\iota}=1\}$, $n_\iota=M_\iota-p_\iota-1$, be a partition of $[0,1]$, called \emph{knot vector}. With this ingredients we are able to define the B-Spline basis $\p{N}_{i,p}$, $i\in\{1,\ldots,M_\iota\}$ on $[0,1]$ via Cox-De Boor's algorithm, cf. \cite{HL:CotrellHughesBazilevs:2009a}. The generalization to $\p{\Omega}$ is realized by considering a tensor product, again denoted by $\p{N}_{i,p}$, where $i=(i_1,\ldots,i_d)$ and $p=(p_1,\ldots,p_d)$ are a multi-indices. For notational simplicity, we define 
	$\iSet:= \{(i_1,\ldots,i_d)\,|\,i_\iota \in \{1,\ldots,M_\iota\}\} $ as the set of multi-indices.
Since the tensor product knot vector $\Xi$ provides a partition of $\p{\Omega}$, it introduces a mesh $\p{\mathcal{Q}}$, and we denote a mesh element by $\p{Q}$, called \emph{cell}.

The B-Spline basis functions parametrize the computational domain  $\Omega$, also called \emph{physical domain}. It is given as image of parameter domain under the \emph{geometrical mapping} $G :\; \p{\Omega} \rightarrow \mathbb{R}^{{d}}$, defined as
\begin{align*}
	    G(\xi) := \sum_{i\in \mathcal{I}} P_i \p{N}_{i,p}(\xi),
\end{align*}
with the control points $P_i \in \mathbb{R}^{{d}}$, $i\in \mathcal{I}$. The image of the mesh $\p{\mathcal{Q}}_h$ under $G$ defines the mesh on $\Omega$, denoted by $\mathcal{Q}_h$ with cells $Q$. Both meshes possess a characteristic mesh size $\p{h}$ and $h$, respectively. More complicated geometries $\Omega$ have to be represented with multiple non-overlapping domains $\Omega\sMP:=G\sMP(\p{\Omega}),k=1,\ldots,N$, called \emph{patches}, where each patch is associated with a different geometrical mapping $G\sMP$. We sometimes call $\overline{\Omega}:=\bigcup_{k=1}^N\overline{\Omega}\sMP$ a \emph{multipatch domain}. Furthermore, we denote the set of all indices $l$  
such that $\Omega^{(k)}$ and $\Omega^{(l)}$ have a common interface $F^{(kl)}$ by ${\mathcal{I}}_{\mathcal{F}}^{(k)}$.
We define the interface $\Gamma\sMP$ of $\Omega\sMP$ as $\Gamma\sMP := \bigcup_{l\in\ifSet\sMP}^N F^{(kl)}$.

We use B-Splines not only for defining the geometry, but also for representing the approximate solution of our PDE. This motivates to define the basis functions in the physical space $\g{N}_{i,p}:=\p{N}_{i,p}~\circ~G^{-1}$ and the corresponding discrete space as
\begin{align}
\label{equ:gVh}
  V_h:=\text{span}\{\g{N}_{i,p}\}_{i\in\iSet}.
\end{align}
Moreover, each function $u_h(x) = \sum_{i\in\mathcal{I}} u_i \g{N}_{i,p}(x)$ 
is associated with the 
 coefficient  vector $\boldsymbol{u} = (u_i)_{i\in\mathcal{I}}$. 
This map is known as \emph{Ritz isomorphism} 
or \emph{IgA isomorphism}  in connection with IgA.
One usually writes this relation as 
$u_h \leftrightarrow \boldsymbol{u}$.
 In the following, we will use the notation $u_h$ for the function and its vector representations.
 If we consider a single patch $\Omega^{(k)}$ of a multipatch domain $\Omega$, we will use the notation $V_{h}^{(k)},\g{N}_{i,p}^{(k)},\p{N}_{i,p}^{(k)}, G^{(k)}, \ldots$ with the analogous definitions. To keep notation simple, we will use $h_k$ and $  \p{h}_k$ instead of $h\sMP$ and $\p{h}\sMP$, respectively.

\section{Galerkin Methods for Isogeometric Analysis}
\label{sec:galerkin-IGA}
In this section we rephrase the variational formulation for the continuous and discontinuous Galerkin method for multipatch IgA systems. 
\subsection{Continuous Galerkin method}
We are considering the finite dimensional subspace 
$V_{h}^{cG}$ of $V_D$, where $V_{h}^{cG}$ is given by
\begin{align*}
	V_h^{cG}:= \{v\,|\, v|_{\Omega\sMP}\in V_h\sMP\} \cap H^1({\Omega}).
\end{align*}
Since, we restrict ourselves to homogeneous Dirichlet conditions, we look for the Galerkin approximate $u_h$ from $V_{D,h}^{cG} \subset V_h^{cG}$, where $V_{D,h}^{cG}$ contains all functions, which vanish on the Dirichlet boundary.
The Galerkin IgA scheme reads as follows:
Find $u_h \in V_{D,h}^{cG}$ such that
  \begin{align}
  \label{HL:equ:ModelDisc}
    a(u_h,v_h) = \left\langle F, v_h \right\rangle \quad \forall v_h \in V_{D,h}.
  \end{align}
There exists a unique IgA solution $u_h \in V_{D,h}^{cG}$ of (\ref{HL:equ:ModelDisc})
that converges to the solution $u \in V_{D}$ of (\ref{equ:ModelVar}) 
if $h$ tends to $0$.
Due to Cea's lemma, the usual discretization error estimates 
in the $H^1$ - norm follow from the corresponding approximation 
error estimates, see \cite{HL:BazilevsVeigaCottrellHughesSangalli:2006a} 
or  \cite{HL:BeiraodaVeigaBuffaSangalliVazquez:2014a}.

\subsection{Discontinuous Galerkin method}
In the dG-IgA scheme, we again use the 
spaces $V_{h}^{(k)}$ of B-Splines defined in \eqref{equ:gVh}, 
whereas now
discontinuities are allowed 
across the patch interfaces $F\sMP[kl]$.
The continuity of the function value and its normal fluxes are then enforced in a weak sense by adding additional terms to the bilinear form.
We define the dG-IgA space
\begin{align}
\label{equ:gVh_glob}
  V_{h}^{dG}:= \{v\,| \,v|_{\Omega^{(k)}}\in V_{h}^{(k)}\},
\end{align}
where $V_{h}^{(k)}$ is defined as in \eqref{equ:gVh}. 
A comprehensive study of dG schemes for FE can be found in \cite{HL:Riviere:2008a} and \cite{HL:PietroErn:2012a}. For an analysis of the dG-IgA scheme, we refer to \cite{HL:LangerToulopoulos:2015a}.
We define $V_{D,h}^{dG}$ as the space of all functions from $V_{h}$ that vanish on the Dirichlet boundary $\Gamma_D$. Having these definitions at hand, we can define the discrete problem based on the Symmetric Interior Penalty (SIP) dG formulation as follows:
Find $u_h \in V_{D,h}^{dG}$ such that
  \begin{align}
  \label{equ:ModelDiscDG}
    a_h(u_h,v_h) = \left\langle F, v_h \right\rangle \quad \forall v_h \in V_{D,h}^{dG},
  \end{align}
where
\begin{align*}
    a_h(u,v) &:= \sum_{k=1}^N a_e^{(k)}(u,v) \quad \text{and} \quad \left\langle F, v \right\rangle:=\sum_{k=1}^N \left( \int_{\Omega^{(k)}}f v^{(k)} dx+\int_{\Gamma_N\sMP} g_N v\sMP \,ds\right), \\
    a^{(k)}_e(u,v) &:=  a^{(k)}(u,v) + s^{(k)}(u,v) + p^{(k)}(u,v),
\end{align*}
and
\begin{align*}
a^{(k)}(u,v) &:= \int_{\Omega^{(k)}}\alpha^{(k)} \nabla u^{(k)} \nabla v ^{(k)} dx,\\
    s^{(k)}(u,v)&:= \sum_{l\in{\mathcal{I}}_{\mathcal{F}}^{(k)}} \int_{F^{(kl)}}\frac{\alpha^{(k)}}{2}\left(\frac{\partial u^{(k)}}{\partial n}(v^{(l)}-v^{(k)})+ \frac{\partial v^{(k)}}{\partial n}(u^{(l)}-u^{(k)})\right)ds,\\ 
    p^{(k)}(u,v)&:= \sum_{l\in{\mathcal{I}}_{\mathcal{F}}^{(k)}} \int_{F^{(kl)}}\frac{\delta \alpha^{(k)}}{h_{kl}}(u^{(l)}-u^{(k)})(v^{(l)}-v^{(k)})\,ds.
\end{align*}
Here the notation $\frac{\partial}{\partial n}$ denotes the derivative 
in the direction of the outer normal vector,
$\delta$ a positive sufficiently large penalty parameter,  and $h_{kl}$ the harmonic average of the adjacent mesh sizes, i.e., $h_{kl}= 2h_k h_l/(h_k + h_l)$.

We equip $V_{D,h}^{dG}$ with the dG-norm 
\begin{align}
\label{HL:dgNorm}
 \left\|u\right\|_{dG}^2  = \sum_{k = 1}^N\left[\alpha^{(k)} \left\|\nabla u^{(k)}\right\|_{L^2(\Omega^{(k)})}^2 + \sum_{l\in{\mathcal{I}}_{\mathcal{F}}^{(k)}} \frac{\delta \alpha^{(k)}}{h_{kl}}\int_{F^{(kl)}} (u^{(k)} - u^{(l)})^2 ds\right].
\end{align}

Furthermore,
we define the bilinear forms
\begin{align*}
 d_h(u,v) = \sum_{k =1}^N d^{(k)}(u,v) \quad \text{where} \quad d^{(k)}(u,v)= a^{(k)}(u,v) + p^{(k)}(u,v),
\end{align*}
for later use.
We note that $\left\|u_h\right\|_{dG}^2 = d_h(u_h,u_h)$.
%
\begin{lemma}
\label{lem:wellPosedDg}
Let $\delta$ be sufficiently large. 
Then there exist two positive constants $\gamma_0$ and $\gamma_1$,
which are independent of $h_k,H_k,\delta,\alpha^{(k)}$ and $u_h$
such that the inequalities
 \begin{align}
 \label{equ:equivPatchNormdG}
  \gamma_0 d^{(k)}(u_h,u_h)\leq a_e^{(k)}(u_h,u_h) \leq \gamma_1 d^{(k)}(u_h,u_h), \quad 
  \forall u_h\in V_{D,h}^{dG}
 \end{align}
are valid for all $k=1,2,\ldots,N$. Furthermore, we have  the inequalities
\begin{align}
 \label{equ:equivNormdG}
 \gamma_0 \left\|u_h\right\|_{dG}^2\leq a_h(u_h,u_h)\leq \gamma_1 \left\|u_h\right\|_{dG}^2, \quad 
  \forall u_h\in V_{D,h}^{dG}.
\end{align}

\end{lemma}
This Lemma is an equivalent statement of Lemma 2.1 in \cite{HL:DryjaGalvisSarkis:2013a} for IgA,
and the proof can be found in \cite{HL:HoferLanger:2016a}. A direct implication of \eqref{equ:equivNormdG} is the well posedness of the discrete problem \eqref{equ:ModelDiscDG} by the Theorem of Lax-Milgram. The consistency of the method together with the interpolation estimates of B-Splines lead to the a-priori error estimate, established in \cite{HL:LangerToulopoulos:2015a}.
%
%
We note that, in \cite{HL:LangerToulopoulos:2015a}, the results were obtained 
for the Incomplete Interior Penalty (IIP) scheme. 
An extension to SIP-dG and the use of  harmonic averages for $h$ and/or $\alpha$ 
are discussed in Remark~3.1 in \cite{HL:LangerToulopoulos:2015a},
see also \cite{HL:LangerMantzaflarisMooreToulopoulos:2015b}.

For both the cG and dG formulation, we choose the B-Spline function $\{\g{N}_{i,p}\}_{i\in\mathcal{I}_0}$ as basis for the space $V_{h}^X,X\in\{cG,dG\}$, where $\mathcal{I}_0$ contains all indices of $\mathcal{I}$, where the corresponding basis functions do not have a support on the Dirichlet boundary. In the cG case, the basis functions on the interface are identified accordingly to obtain a conforming subspace of $V_{D}$. For the remainder of this paper, we drop the superscript $X\in\{cG,dG\}$ and use the symbol $V_h$ for both formulations. Depending on the considered formulation, one needs to use the right space $V_{h}^X,X\in\{cG,dG\}$. The IgA schemes \eqref{HL:equ:ModelDisc} and \eqref{equ:ModelDiscDG} are equivalent to the system of linear IgA equations
\begin {align}
\label{equ:Ku=f_DG}
  \boldsymbol{K} \boldsymbol{u} = \boldsymbol{f},
\end{align}
where 
$\boldsymbol{K} = (\boldsymbol{K}_{i,j})_{i,j\in {\mathcal{I}}_0}$
,
$\boldsymbol{f}= (\boldsymbol{f}_i)_{i\in {\mathcal{I}}_0}$
denote the stiffness matrix and the load vector, respectively,
with 
$ \boldsymbol{K}_{i,j} = a(\g{N}_{j,p},\g{N}_{i,p})$ or $ \boldsymbol{K}_{i,j} = a_h(\g{N}_{j,p},\g{N}_{i,p})$
and 
$\boldsymbol{f}_i = \left\langle F, \g{N}_{i,p} \right\rangle$, and
$\boldsymbol{u}$ is the vector representation of $u_h$. 

\section{IETI-DP methods and their implementation}
\label{sec:IETI-method}
In this section, we recall the main ingredients for the cG-IETI-DP and dG-IETI-DP method. We focus mainly on the implementation, since this is the relevant part for parallelization.
\subsection{Derivation of the method}
\label{sec:derivation}

A rigorous and formal definition of the cG-IETI-DP and dG-IETI-DP method is quite technical and not necessary for the parallelization, which is the purpose of this paper. Therefore, we are not going to present the whole derivation of each method. We will give a general description, which is valid for both methods. For a detailed derivation, we refer to \cite{HL:HoferLanger:2016b} and \cite{HL:HoferLanger:2016a}. 

The first step is to introduce additional dofs on the interface to decouple the local problems and incorporate their connection via Lagrange multipliers $\boldsymbol{\lambda}$. This is quite straightforward in the case of the cG formulation, but more involved in the dG case. In any of the two cases, we can equivalently rewrite \eqref{equ:Ku=f_DG} as: Find $(u,\boldsymbol{\lambda}) \in V_{h,e} \times U$ such that
    \begin{align}
    \label{equ:saddlePointSing}
     \MatTwo{K_e}{B^T}{B}{0} \VecTwo{u}{\boldsymbol{\lambda}} = \VecTwo{f}{0},
    \end{align}
where $V_{h,e}\supset V_{h}$, is the decoupled space with additional dofs and $U$ is the set of Lagrange multipliers. The jump operator $B$ enforces the ``continuity'' of the solution $u$ in the sense that $\ker(B) \equiv  V_{h}$. The matrix $K_e$ is the block diagonal matrix of the patch local stiffness matrices $K\sMP$, i.e., $K_e = \text{diag}(K\sMP)$. Since $B$ only acts on the patch interface dofs, we first can reorder the stiffness matrix in the following way
  \begin{align*}
    K\sMP = \MatTwo{K_{BB}\sMP}{K_{BI}\sMP}{K_{IB}\sMP}{K_{II}\sMP}, \qquad f\sMP = \VecTwo{f_B\sMP}{f_I\sMP}
  \end{align*}
and then consider only the Schur complement representation:
Find $(u_B,\boldsymbol{\lambda}) \in W \times U$ such that
    \begin{align}
    \label{equ:saddlePointSingSchur}
     \MatTwo{S_e}{B_B^T}{B_B}{0} \VecTwo{u_B}{\boldsymbol{\lambda}} = \VecTwo{g}{0},
    \end{align}
where $S_e=\text{diag}(S_e\sMP)$,  $S_e\sMP= K_{BB}\sMP- K_{BI}\sMP(K_{II}\sMP)^{-1}K_{IB}\sMP$  and $g\sMP= f_B -  K_{BI}\sMP (K_{II}\sMP)^{-1}f_I\sMP$. The space $W$ is the restriction of $V_{h,e}$ to the interface. For completeness, we denote its ``continuous'' representation as $\widehat{W}$, i.e., $\ker(B_B) = \widehat{W}$. Equation \eqref{equ:saddlePointSingSchur} is valid for the cG-IETI-DP and dG-IETI-DP method, but the matrix $K_e$ has difference entries and the number of boundary dofs (subscript $B$) is different. Fortunately, this does not change the way how the algorithm is implemented and parallelized. In the following, we will drop the subscript $B$ in $u_B$ and $B_B$ for notational simplicity.

The matrix $S_e$ is not invertible and, hence, we cannot build the Schur complement system of \eqref{equ:saddlePointSingSchur}. To overcome this, we introduce an intermediate space $\widetilde{W}$, such that $\widehat{W} \subset \widetilde{W} \subset W$, and $S_e$ restricted to $\widetilde{W}$, denoted by $\widetilde{S}$, is invertible. We introduce primal variables as a set $\Psi\subset \widehat{W}^*$ and define the spaces
\begin{equation*}
\widetilde{W} := \{w\in W:  \psi(w^{(k)}) = \psi(w^{(l)}), \forall\psi \in \Psi, \forall k>l  \}
\end{equation*}
and
\begin{equation*}
W_{\Delta} := \prod_{k=1}^N W_{\Delta}^{(k)},\text{ with} \quad W_{\Delta}^{(k)}:=\{w^{(k)}\in W^{(k)}:\, 
\psi(w^{(k)}) =0\; \forall\psi \in \Psi\}.
\end{equation*}
Moreover, we introduce the space $W_{\Pi} \subset \widehat{W}$ such that
$\widetilde{W} = W_{\Pi} \oplus W_{\Delta}.$
We call $W_{\Pi}$ \emph{primal space} and $W_{\Delta}$ \emph{dual space}.  Typically, the set $\Psi$ corresponds to ``continuous'' vertex values, edge averages and/or face averages. 

Since $\widetilde{W} \subset W$, there is a natural embedding $\widetilde{I}: \widetilde{W} \to W$. 
Let the jump operator restricted to $\widetilde{W}$ be
    $\widetilde{B} := B\widetilde{I} : \widetilde{W} \to U^*.$
Now we are in the position to reformulate problem \eqref{equ:saddlePointSingSchur} in the space $\widetilde{W}$ as follows:
Find $(u,\boldsymbol{\lambda}) \in \widetilde{W} \times U:$
    \begin{align}
    \label{equ:saddlePointReg}
     \MatTwo{\widetilde{S}}{\widetilde{B}^T}{\widetilde{B}}{0} \VecTwo{u}{\boldsymbol{\lambda}} = \VecTwo{\widetilde{g}}{0},
    \end{align}
    where $\widetilde{g} := \widetilde{I}^T g$, and $\widetilde{B}^T= \widetilde{I}^T B^T$. Here, $\widetilde{I}^T: W^* \to \widetilde{W}^*$ denotes the adjoint of $\widetilde{I}$.

By construction, $\widetilde{S}$ is SPD on $\widetilde{W}$. Hence, we can define the Schur complement $F$ 
and the corresponding right-hand side
as follows:
\begin{align*}
    F:= \widetilde{B} \widetilde{S}^{-1}\widetilde{B}^T, \quad d:= \widetilde{B}\widetilde{S}^{-1} \widetilde{g}.
\end{align*}
Hence, the saddle point system \eqref{equ:saddlePointReg} is equivalent to 
the Schur complement problem:
\begin{align}
   \label{equ:SchurFinal}
      \text{Find } \boldsymbol{\lambda} \in U: \quad F\boldsymbol{\lambda} = d.
\end{align}
Equation \eqref{equ:SchurFinal} is solved by means of the PCG algorithm, but it requires an appropriate preconditioner in order to obtain an efficient solver.

Recalling the definition of $S_e = \text{diag}(S_e^{(k)})_{k=1}^N$, we define the scaled Dirichlet preconditioner 
	$M_{sD}^{-1} := B_D S_e B_D^T,$
where $B_D$ is a scaled version of the jump operator $B$. 
The scaled jump operator
$B_D$ is defined such that the operator enforces the constraints
 \begin{align*}
  {\delta^\dagger}^{(l)}_j(\boldsymbol{u}^{(k)})^{(k)}_i -  {\delta^\dagger}^{(k)}_i(\boldsymbol{u}^{(l)})^{(k)}_j = 0 \quad\forall (i,j)\in B_e(k,l), \;\forall l\in{\mathcal{I}}_{\mathcal{F}}^{(k)},
 \end{align*}
and
 \begin{align*}
  {\delta^\dagger}^{(l)}_j(\boldsymbol{u}^{(k)})^{(l)}_i -  {\delta^\dagger}^{(k)}_i(\boldsymbol{u}^{(l)})^{(l)}_j = 0 \quad\forall (i,j)\in B_e(l,k), \;\forall l\in{\mathcal{I}}_{\mathcal{F}}^{(k)},
 \end{align*}
 where, for $(i,j)\in B_e(k,l)$,
       $ {\delta^\dagger}^{(k)}_i:= \rho^{(k)}_i/\sum_{l\in{\mathcal{I}}_{\mathcal{F}}^{(k)}} \rho^{(l)}_j $
 is an appropriate scaling. One can show, that  the preconditioned system has a quasi-optimal condition number bound with respect to $H/h:=\max_k(H_k/h_k)$, i.e.,
 \begin{align}
 \label{equ:kappa}
 \kappa(M_{sD}^{-1}F_{|\text{ker}(\widetilde{B}^T)}) \leq C (1+\log(H/h))^2,
 \end{align}
 for both versions, see \cite{HL:HoferLanger:2016b}, \cite{HL:Hofer:2016a} and \cite{HL:VeigaChoPavarinoScacchi:2013a}. Moreover, numerical examples show also robustness with respect to jumps in the diffusion coefficient and only a weak dependence on the B-Spline degree $p$, see, e.g., \cite{HL:HoferLanger:2016a}, \cite{HL:HoferLanger:2016b} and \cite{HL:VeigaChoPavarinoScacchi:2013a}.

\subsection{Implementation of the algorithm}
\label{sec:implementation}
  Since $F$ is symmetric and  positive definite on $\widetilde{U}$, we can solve
   the linear system $F \mathbf{\boldsymbol{\lambda}} = d$ by means of the PCG algorithm, 
where we use $M_{sD}^{-1}$ as preconditioner. The PCG does not require an explicit representation of the matrices $F$ and $M_ {sD}^{-1}$, since we just need their application to a vector.
There are different ways to provide an efficient implementation. We will follow the concept of the energy minimizing primal subspaces. The idea is to split the space $\widetilde{W}$ into $\widetilde{W}_{\Pi} \oplus \prod \widetilde{W}_{\Delta}^{(k)}$, such that $\widetilde{W}_{\Delta}^{(k)} {\perp_S} \widetilde{W}_{\Pi}$ for all $k$, i.e., we choose $\widetilde{W}_{\Pi}:=\widetilde{W}_{\Delta}^{\perp_S}$, see, e.g., \cite{HL:Pechstein:2013a} and \cite{HL:Dohrmann:2003a}.
By means of this choice, the operators $\tilde{S}$ and $\tilde{S}^{-1}$ have the following forms
   \begin{align*}
    \widetilde{S} &= \MatTwo{S_{\Pi \Pi}}{0}{0}{S_{\Delta \Delta}} \text{ and }  \widetilde{S}^{-1} = \MatTwo{S_{\Pi \Pi}^{-1} }{0}{0}{S_{\Delta \Delta}^{-1}},
    \end{align*}
 where $S_{\Pi \Pi}$ and $S_{\Delta \Delta}$ are the restrictions of $\widetilde{S}$ to the corresponding subspaces. We note that $S_{\Delta \Delta}$ can be seen as a block diagonal operator, i.e., $S_{\Delta \Delta} = \text{diag}(S_{\Delta \Delta}^{(k)})$. 
 
 The application of $F$ and $M_{sD}^{-1}$ is summarized in Algorithm~\ref{HL:alg:applyF}.


\subsubsection{Constructing a basis for the primal subspace}
  First we need to provide an appropriate  local basis $\{\widetilde{\phi}_j\}_j^{n_{\Pi}}$ for $\widetilde{W}_{\Pi}$, where $n_{\Pi}$  
  is the number of primal variables. We request from the basis that it has to be nodal with respect to the primal variables, i.e.,
  $\psi_i(\widetilde{\phi}_j) = \delta_{i,j},$ for $i,j \in\{1,\ldots,n_{\Pi}\}.$ In order to construct such a basis, we introduce the constraint matrix $C^{(k)}: W^{(k)}\to \mathbb{R}^{n_{\Pi}^{(k)}}$ for each patch $\Omega^{(k)}$ which realizes the primal variables, i. e., 
	  $  (C^{(k)} v)_j = \psi_{i(k,j)}(v)$ for $ v\in W$ and $j\in\{1,\ldots,n_{\Pi}^{(k)}\},$
   where $n_{\Pi}^{(k)}$ is the number of primal variables associated with $\Omega^{(k)}$ and $i(k,j)$ the global index of the $j$-th primal variable on $\Omega^{(k)}$.
   For each patch $k$, the basis functions $\{\widetilde{\phi}_j^{(k)}\}_{j=1}^{n_{\Pi}^{(k)}}$ of $\widetilde{W}_{\Pi}^{(k)}$ are the solution of the system
\begin{align}
     \label{HL:equ:KC_basis}
     \MatThree{K_{BB}^{(k)}}{K_{BI}^{(k)}}{{C^{(k)}}^T}{K_{IB}^{(k)}}{K_{II}^{(k)}}{0}{C^{(k)}}{0}{0}\VecThree{\widetilde{\phi}_j^{(k)}}{\cdot}{\widetilde{\mathbf{\mu}}_j^{(k)}} = \VecThree{0}{0}{\mathbf{e}_j^{(k)}},
\end{align}
where  $\mathbf{e}_j^{(k)} \in \mathbb{R}^{n_{\Pi}^{(k)}}$ is the $j$-th unit vector. Here we use an equivalent formulation with the system matrix $K\sMP$
For each patch $k$, the LU factorization of this matrix is computed and stored.


\paragraph{Application of ${S_{\Delta \Delta}^{(k)}}^{-1}:$ } 
     The application of ${S_{\Delta \Delta}^{(k)}}^{-1}$ corresponds to solving a local Neumann problem in the space $\widetilde{W}_{\Delta}$, i.e.,
       $S^{(k)} w^{(k)}= f_{\Delta}^{(k)}$ with the constraint $C^{(k)} w^{(k)} = 0$.
          This problem can be rewritten as a saddle point problem in the form
     \begin{align*}
     \label{HL:equ:SC_solution}
      \MatThree{K_{BB}^{(k)}}{K_{BI}^{(k)}}{{C^{(k)}}^T}{K_{IB}^{(k)}}{K_{II}^{(k)}}{0}{C^{(k)}}{0}{0}\VecThree{w^{(k)}}{\cdot}{\cdot} = \VecThree{f_{\Delta}^{(k)}}{0}{0}.
     \end{align*}
     From (\ref{HL:equ:KC_basis}), the LU factorization of the matrix is already available.
  \paragraph{Application of ${\mathbf{S}_{\Pi\Pi}^{(k)}}^{-1}:$}
     The matrix $\mathbf{S}_{\Pi \Pi}$ can be assembled from the patch local matrices $\mathbf{S}_{\Pi\Pi}^{(k)}$.
     Let $\{\widetilde{\phi}_j^{(k)}\}_{j=1}^{n_{\Pi}^{(k)}}$ be the basis of $\widetilde{W}_{\Pi}^{(k)}$. 
     The construction of $\{\widetilde{\phi}_j^{(k)}\}_{j=1}^{n_{\Pi}^{(k)}}$ in (\ref{HL:equ:KC_basis}) provides 
     \begin{align*}
      \left(\mathbf{S}_{\Pi\Pi}^{(k)}\right)_{i,j} &= \left\langle S^{(k)} \widetilde{\phi}_i^{(k)}, \widetilde{\phi}_j^{(k)} \right\rangle = -\left\langle{C^{(k)}}^T \widetilde{\mathbf{\mu}}_i^{(k)},\widetilde{\phi}_j^{(k)}\right\rangle = -\left\langle \widetilde{\mathbf{\mu}}_i^{(k)}, C^{(k)} \widetilde{\phi}_j^{(k)} \right\rangle\\
	& = -\left\langle \widetilde{\mathbf{\mu}}_i^{(k)}, \mathbf{e}_{j} \right\rangle^{(k)}= - \left(\widetilde{\mathbf{\mu}}_i^{(k)}\right)_j,
     \end{align*}
     where $i,j\in\{1,\ldots,n_{\Pi}^{(k)}\}$. Therefore, we can reuse the Lagrange multipliers $\widetilde{\mathbf{\mu}}_i^{(k)}$ obtained in (\ref{HL:equ:KC_basis}), and can assemble $\mathbf{S}_{\Pi \Pi}^{(k)}$ from them. Once $\mathbf{S}_{\Pi \Pi}$ is assembled, the  LU factorization can be calculated and stored.

\subsubsection{Application of $\widetilde{I}$ and $\widetilde{I}^T$}
The last building block is the embedding $\widetilde{I}: \widetilde{W}\to W$ and its adjoint $\widetilde{I}^T: W^*\to \widetilde{W}^*$. Recall the direct splitting $W\sMP = W_\Delta\sMP \oplus W_{\Pi}\sMP$. Let us denote by $\Phi\sMP=[\widetilde{\phi}_1\sMP,\ldots,\widetilde{\phi}_{n_\Pi\sMP}\sMP]$ the coefficient representation of the basis for $W_{\Pi}\sMP$. Given the primal part $\boldsymbol{w}_\Pi$ of a function in $\widetilde{W}$, we obtain its restriction to $\widetilde{W}_{\Pi}\sMP$ via an appropriately defined restriction matrix $\boldsymbol{R}\sMP$, i.e. $\boldsymbol{w}_\Pi\sMP = \boldsymbol{R}\sMP\boldsymbol{w}_\Pi$. The corresponding function is then given by $w_\Pi\sMP = \Phi\sMP\boldsymbol{R}\sMP\boldsymbol{w}_\Pi\sMP$.

Following the lines in \cite{HL:Pechstein:2013a}, we can formulate the operator $\widetilde{I}: \widetilde{W}\to W$ as
\begin{align*}
	\begin{bmatrix}
		\boldsymbol{w}_\Pi\\
		w_\Delta
	\end{bmatrix} 
	\mapsto w:= \Phi \boldsymbol{R} \boldsymbol{w}_\Pi + w_\Delta,
\end{align*}
where $\Phi$ and $\boldsymbol{R}$ are block versions of $\Phi\sMP$ and $\boldsymbol{R}\sMP$, respectively. The second function is its adjoint operation $\widetilde{I}^T: W^*\to \widetilde{W}^*$. It can be realized in the following way
\begin{align*}
f \mapsto
	\begin{bmatrix}
		\boldsymbol{f}_\Pi\\
		f_\Delta
	\end{bmatrix} 
	=
	\begin{bmatrix}
		\boldsymbol{A} \Phi^T \boldsymbol{f}\\
		f - C^T \Phi^T f
	\end{bmatrix}, 
\end{align*}
where $\boldsymbol{A}$ is the corresponding assembling operator to $\boldsymbol{R}$, i.e., $\boldsymbol{A} = \boldsymbol{R}^T$.   A more extensive discussion and derivation can be found in \cite{HL:Pechstein:2013a}.

   \begin{algorithm}
   \caption{Algorithm for calculating $\boldsymbol{\nu} = F\mathbf{\boldsymbol{\lambda}}$  and 
    $\mathbf{\boldsymbol{\nu}} = M_{sD}^{-1}\mathbf{\boldsymbol{\lambda}}$ for given $\mathbf{\boldsymbol{\lambda}} \in U$}
   \label{HL:alg:applyF}
   \begin{algorithmic}
   \algblock{Begin}{End}
   \Procedure{$F$}{$\mathbf{\boldsymbol{\lambda}}$}
      \State Application of $B^T:$ $\{f^{(k)}\}_{k=1}^ N = B^T\mathbf{\boldsymbol{\lambda}}$
      \State Application of $\widetilde{I}^T:$ $\{\mathbf{f}_{\Pi},\{f_{\Delta}^{(k)}\}_{k=1}^ N\} = \widetilde{I}^T\left(\{f^{(k)}\}_{k=1}^ N\right)$
      \State  Application of $\widetilde S^{-1}:$
      \Begin
	    \State  $\mathbf{w}_{\Pi} = \mathbf{S}_{\Pi \Pi}^{-1} \mathbf{f}_{\Pi}$
	    \State  $w_{\Delta}^{(k)} = {S_{\Delta \Delta}^{(k)}}^{-1}f_{\Delta}^{(k)} \quad \forall k=1,\ldots, N$ 
      \End
      \State  Application of $\widetilde{I}:$ $\{w^{(k)}\}_{k=1}^ N = \widetilde{I}\left(\{\mathbf{w}_{\Pi},\{w_{\Delta}^{(k)}\}_{k=1}^ N\} \right)$
      \State  Application of $B:$ $ \boldsymbol{\nu} = B\left( \{w^{(k)}\}_{k=1}^ N \right)$
   \EndProcedure
   \Procedure{$M_{sD}^{-1}$}{$\mathbf{\boldsymbol{\lambda}}$}
      \State Application of $B_D^T:$ $\{w^{(k)}\}_{k=1}^N = B_D^T\mathbf{\boldsymbol{\lambda}}$
      \State  Application of $S_e:$
      \Begin
	    \State Solve $K^{(k)}_{II} x^{(k)} = -K^{(k)}_{IB}w^{(k)}$ $\quad \forall k=1,\ldots, N$ 
	    \State $v^{(k)} = K^{(k)}_{BB} w^{(k)} + K^{(k)}_{BI}x^{(k)}$. $\quad \forall k=1,\ldots, N$ 
      \End
      \State  Application of $B_D:$ $ \mathbf{\boldsymbol{\nu}} = B_D\left( \{v^{(k)}\}_{k=1}^N \right)$
   \EndProcedure
  \end{algorithmic}
  \end{algorithm}
\section{Parallelization of the building blocks}
\label{sec:para}
Here we investigate how the single operations can be executed in parallel in a distributed memory setting. The parallelization of the method is performed with respect to the patches, i.e., one or several patches are assigned to a processor. The required communication has to be understood as communication between patches, which are assigned to different processors. The majority of the used MPI methods are performed in its non-blocking version. We aim at overlapping computations with communications wherever possible. 

\subsection{Parallel version of PCG}
We solve $F \boldsymbol{\lambda} = d$ with the preconditioned CG method. This requires a parallel implementation of CG, where we follow the approach presented in Section 2.2.5.5 in \cite{HL:Pechstein:2013a}, see also \cite{HL:DouglasHaaseLanger:2003a}. This approach is based on the concept of accumulated and distributed vectors. We say a vector $\boldsymbol{\lambda}_{acc}=[\boldsymbol{\lambda}_{acc}\sMP[q]]$ is an \emph{accumulated} representation of $\boldsymbol{\lambda}$, if $\boldsymbol{\lambda}_{acc}\sMP[q](k_q(i)) = \boldsymbol{\lambda}(i)$, where $i$ is the global index corresponding to the  local index $k_q(i)$ on processor $q$.
On the contrary, $\boldsymbol{\lambda}_{dist}=[\boldsymbol{\lambda}_{dist}\sMP[q]]$ is a \emph{distributed} representation of $\boldsymbol{\lambda}$, if the sum of all processor local contributions give the global vector, i.e.,
	$\boldsymbol{\lambda}_{dist}(i) = \sum_{q} \boldsymbol{\lambda}_{dist}\sMP[q](k_q(i))$.
Hence, each processor only holds the part of $\boldsymbol{\lambda}$, which belongs to its patches, either in a distributed or accumulated description. The Lagrange multipliers and the search direction of the CG are represented in the accumulated setting, whereas the residual is given in the distributed representation. In order to achieve the accumulated representation, information exchange between the neighbours of a patch is required. This is done after applying the matrix and the preconditioner, respectively and implemented via \verb+MPI_Send+ and \verb+MPI_Recv+ operations.

The last aspect in the parallel CG implementation is the realization of scalar products. Given a distributed representation $u_{dist}$ of $u$ and an accumulated representation of $v_{acc}$ of $v$, the scalar product $(u,v)_{l^2}$ is then given by
	$(u,v)_{l^2}= \sum_{q} (u_{dist}\sMP[q],v_{acc}\sMP[q])_{l^2}$,
i.e., first the local scalar products are formed, globally added,  and distributed with \verb+MPI_Allreduce+  .

\subsection{Assembling}
The assembling routine of the IETI-DP algorithm consists of the following steps:
\begin{enumerate}
	\item Assemble the patch local stiffness matrices and right hand side,
	\item assemble the system matrix in \eqref{HL:equ:KC_basis} and calculating its LU-factorization,
	\item assemble $S_{\Pi\Pi}$ and calculating its LU-factorization,
	\item calculate the LU-factorization of $K_{II}\sMP$ ,
	\item calculate the right hand side $\{g_\Pi, g_\Delta\} =\widetilde{g}\in\widetilde{W}^*$, with $g\sMP=f_B -  K_{BI}\sMP (K_{II}\sMP)^{-1}f_I\sMP$.
\end{enumerate}
Most of the tasks are completely independent of each other and, hence, can be performed in parallel. Only the calculation of $S_{\Pi\Pi}$ and $\widetilde{g} = \widetilde{I}^T g$ require communication, which will be handled in Section~\ref{sec:para_AccDist}. 

The LU-factorization of $S_{\Pi\Pi}$ is only required at one processor, since it has to be solved only once per CG iteration. According to \cite{HL:KlawonnLanserRheinbachStengelWellein:2015a}, it is advantageous to distribute this matrix to all other processors in order to reduce communication in the solver part, see  \cite{HL:KlawonnRheinbach:2010a} and references therein for improving scalability based on a different approach.
In the current paper, we investigate cases, where one, several and all processor hold the LU-factorization of $S_{\Pi\Pi}$.  Therefore, each processor is assigned to exactly one holder of $S_{\Pi\Pi}$. This relation is implemented by means of an additional MPI communicator. 

We note that, for extremely large scale problems with $\geq 10^5$ subdomains, one has to consider different strategies dealing with $S_{\Pi\Pi}^{-1}$. Most commonly one uses AMG and solves $S_{\Pi\Pi} u_\Pi = f_\Pi$ in an inexact way, see, e.g., \cite{HL:KlawonnLanserRheinbach:2015a} and \cite{HL:KlawonnRheinbach:2007a}. When considering a moderate number of patches, i.e., $10^3 - 10^4$,  the approach using the LU-factorization of $S_{\Pi\Pi}$ is the most efficient one. In this paper, we restrict ourselves to this case.

The patch local matrix $S_{\Pi\Pi}\sMP$ is obtained as a part of the solutions of \eqref{HL:equ:KC_basis} and the assembling of the global matrix $S_{\Pi\Pi}$ is basically a \verb+MPI_gatherv+ operation. In the case where all processors hold $S_{\Pi\Pi}$ we use \verb+MPI_allgatherv+. If several processors hold the LU factorization, we just call \verb+MPI_gatherv+ on each of these processors. A different possibility would be to first assemble $S_{\Pi\Pi}$  on one patch, distribute it to the other holders and then calculate the LU-factorization on each of the processors. 
\subsection{Solver and Preconditioner}
\label{sec:para_AccDist}
More communication is involved in the solver part. According to Algorithm\,\ref{HL:alg:applyF}, we have to perform the following operations:
\begin{enumerate}
	\item application of $B$ and $B^T$ and its scaled versions
	\item application of $\widetilde{I}$ and $\widetilde{I}^T$
	\item application of $\widetilde{S}^{-1}$
	\item application of $S^{-1}$
\end{enumerate}
The only operations which require communication are $\widetilde{I}$ and $\widetilde{I}^T$. To be more precise, the communication is hidden in the operators $\boldsymbol{A}$ and $\boldsymbol{R}$, see Section~\ref{sec:implementation}, all other operations are block operations, where the corresponding matrices are stored locally on each processor. In principle, their implementation is given by accumulating and distributing values.  The actual implementation depends on how many processors hold the coarse grid problem.

\begin{figure}[h]
  \begin{subfigure}{0.48\textwidth}
        \center{\includegraphics[width=0.5\textwidth]{\PathToPic 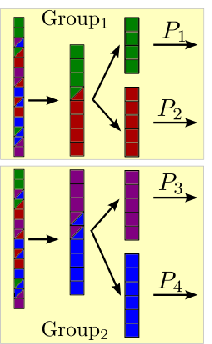}}
  	\caption{Distribution operation}
  \end{subfigure}
    \begin{subfigure}{0.48\textwidth}
        \center{\includegraphics[width=0.5\textwidth]{\PathToPic 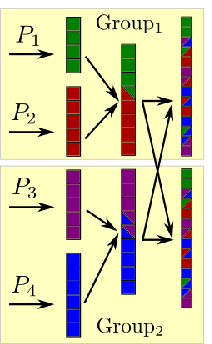}}
  	\caption{Assembling operation}
  \end{subfigure}
  \caption{Distribution and  assembling operation, illustrated for four processors, partitioned into two groups corresponding to two $S_{\Pi\Pi}^{-1}$ holder.}
  \label{fig:DistAcc}
\end{figure}

In order to implement $\widetilde{I}$, we need the distribution operation $\boldsymbol{R}$. If all processors hold $S_{\Pi\Pi}$, this operations reduces to just extracting the right entries. Hence it is local and no communication is required. Otherwise, each holder of $S_{\Pi\Pi}$ reorders and duplicates the entries of $\boldsymbol{w}_\Pi$ in such a way, that all entries corresponding to the patches of a single slave are in a contiguous block of memory. Then we utilize the \verb+MPI_scatter+ method to distribute only the necessary data to all slave processors. See Figure~\ref{fig:DistAcc} for an illustration. 

We arrive at the implementation of $\widetilde{I}^T$. Each processor stores the values of $\boldsymbol{w}\sMP_{\Pi}$ in a vector $\widetilde{\boldsymbol{w}}\sMP_{\Pi}$ of length $n_\Pi$ already in such a way, that 
	$\sum_{k=1}^N \widetilde{\boldsymbol{w}}\sMP_{\Pi} = \boldsymbol{w}_{\Pi}.$
Storing the entries in this way enables the use the MPI reduction operations to efficiently assemble the local contributions. If only one processor holds the coarse problem, we use the \verb+MPI_Reduce+ method to perform this operation. Similarly, if all processors hold $S_{\Pi\Pi}$, we utilize the \verb+MPI_Allreduce+ method. If several processors have the coarse grid problem, we use a two level approach. First, each master processor collects the local contributions from its slaves using the \verb+MPI_Reduce+ operation. In the second step, all the master processors perform an \verb+MPI_Allreduce+ operation to accumulate the contributions from each group and simultaneously distribute the result under them. This procedure is visualized in Figure~\ref{fig:DistAcc}.

\section{Numerical examples}
\label{sec:num_ex}
We consider the model problem \eqref{equ:ModelStrong} in the two dimensional computational domain $\Omega = (0,1)^2$ formed by $32\times32 = 1024$ patches. Each of them is a square arranged in a uniform grid. For the three dimensional case we consider the domain $\Omega=(0,1)^2\times(0,2)$, partitioned into $8\times8\times16$ regular cubes. Note that, in IgA framework, we cannot choose the number of subdomains as freely as in the finite element case since they are fixed by the geometry. Therefore, the number of 1024 subdomains stays constant throughout the tests. Since we are interested in the parallel scalability of the proposed algorithms, we assume for simplicity homogeneous diffusion coefficients $\alpha \equiv 1$. In all tests we consider the smooth right hand side $f(x,y) = 20\pi^2\sin(4\pi(x+0.4))\sin(2\pi(y+0.3))$, corresponding to the exact solution $u(x,y) = \sin(4\pi(x+0.4))\sin(2\pi(y+0.3))+x+y$. For the discretization, we use tensor  B-Spline spaces $V_h$ of different degree $p$. We increase the  B-Spline degree in such a way that the number of knots stay the same, i.e., the smoothness of $V_h$ increases. 

We investigate the scaling behaviour of the cG-IETI-DP and dG-IETI-DP method. Although, we consider also the dG variant, we restrict ourselves to matching meshes. Otherwise, it would not be possible to compare the two methods. Moreover, some patches would have a significant larger number of dofs, which leads to load imbalances and affects the scaling in a negative way.  The domain is refined in a uniform way by inserting a single knot for each dimension on each knot span. We denote by $H_k$ the patch diameter and by $h_k$ the characteristic mesh size on $\Omega\sMP$.  The set of  primal variables is chosen by continuous patch vertices and interface averages for the two dimensional setting. For the three dimensional examples, we choose only continuous edge averages in order to keep the number of primal variables small.

The preconditioned conjugate gradient method is used to solve \eqref{equ:SchurFinal} with the scaled Dirichlet preconditioner $M_{sD}^{-1}$. We choose zero initial guess and a relative reduction of the residual of $10^{-8}$. For solving the local systems and the coarse grid problem, a direct solver is used.

The algorithm is realized in the isogeometric open source C++ library G+SMO \cite{gismoweb}, which is based on the Eigen library \cite{HL:eigenweb}. We utilize the PARDISO 5.0.0 Solver Project \cite{HL:PARDISO500} for performing the LU factorizations. The code is compiled with the \verb+gcc 4.8.3+ compiler with optimization flag \verb+-O3+. For the communication between the processors, we use the \verb+MPI 2+ standard with the \verb+OpenMPI 1.10.2+ implementation. The results are obtain on the RADON1 cluster at Linz. We use 64 out of 68 available nodes, each equipped with 2x Xeon E5-2630v3 ``Haswell'' CPU (8 Cores, 2.4Ghz, 20MB Cache) and 128 GB RAM. This gives the total number of 1024 available cores. 

We investigate two quantities, the assembling phase and the solving phase. In the assembling phase, we account for the time used for 
\begin{itemize}
	\item assembling the local matrices and right hand sides,
	\item LU-factorization of $K_{II}$,
	\item LU-factorization of $\MatTwo{K}{C^T}{C}{0}$,
	\item calculation of $\widetilde{\Phi}$ and $\widetilde{\mu}$,
	\item assembling the coarse grid matrix $S_{\Pi\Pi}$ and calculation of its LU factorization.
\end{itemize}
As already indicated in Section~\ref{sec:para}, $S_{\Pi\Pi}$ is only assembled on certain processors. The solving phase consists of the CG algorithm for \eqref{equ:SchurFinal} and the back-substitution to obtain the solution from the Lagrange multipliers. The main ingredients are
\begin{itemize}
	\item application of $F$,
	\item application of $M_{sD}^{-1}$.
\end{itemize}
In Section~\ref{sec:weak} and Section~\ref{sec:strong}, we study the weak and strong scaling behaviour for the cG-IETI-DP and the dG-IETI-DP method. In this two sections, we assume that only one processor holds the coarse grid matrix $S_{\Pi\Pi}$. The comparison of having a different number of $S_{\Pi\Pi}$ holders is done in Section~\ref{sec:diffHolder}.

\subsection{Weak scaling}
\label{sec:weak}
In this  subsection we investigate the weak scaling behaviour, i.e., the relation of problem size and number of processors is constant. In each refinement step we multiply the number of used cores by $2^{d}, d\in\{2,3\}$. The ideal behaviour  would be a constant time for each refinement.

\begin{figure}[h]
  \begin{subfigure}{0.328\textwidth}
        \includegraphics[width=0.99\textwidth]{\PathToPic 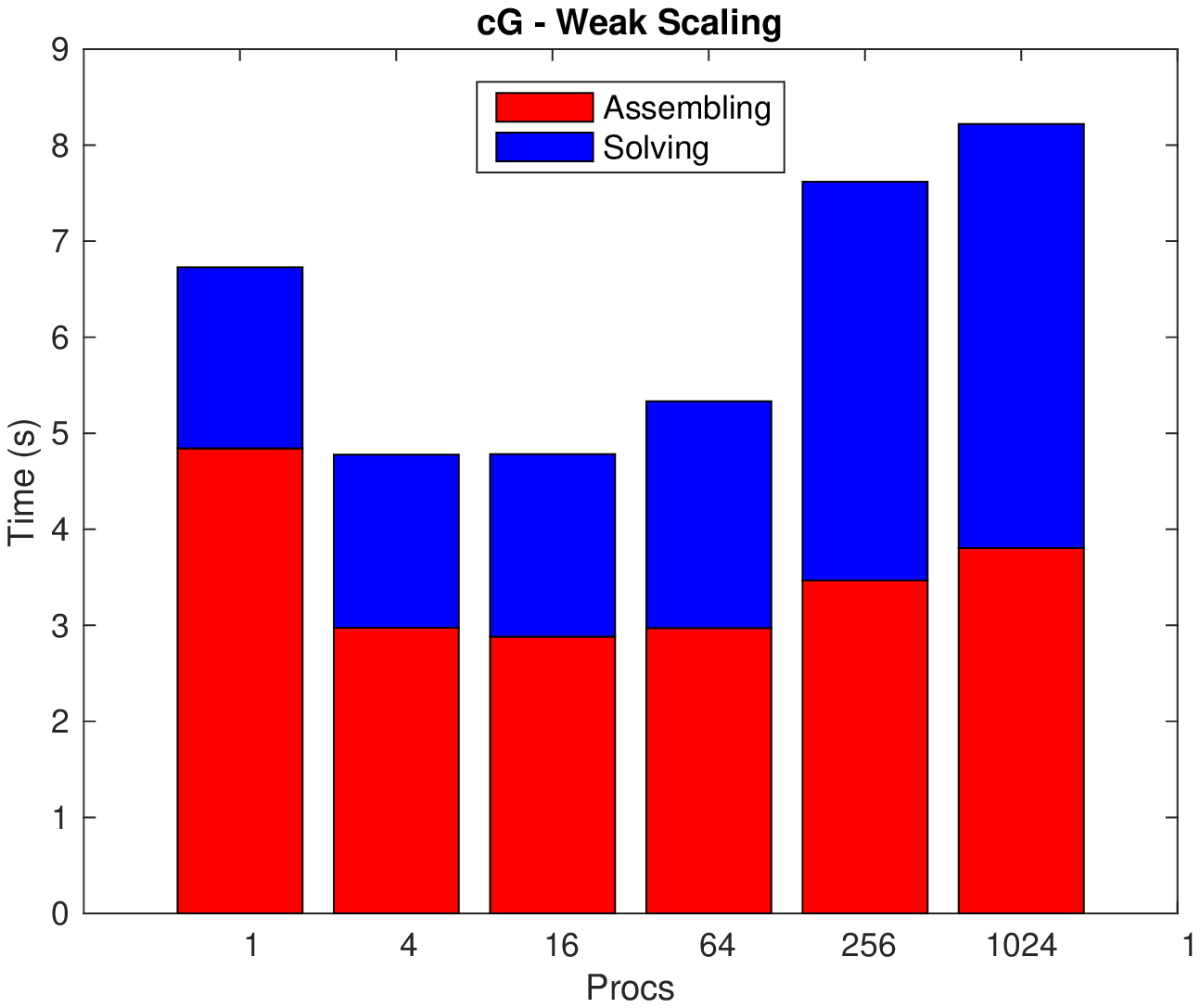}
  	\caption{$p=2$}
  \end{subfigure}
    \begin{subfigure}{0.328\textwidth}
        \includegraphics[width=0.99\textwidth]{\PathToPic 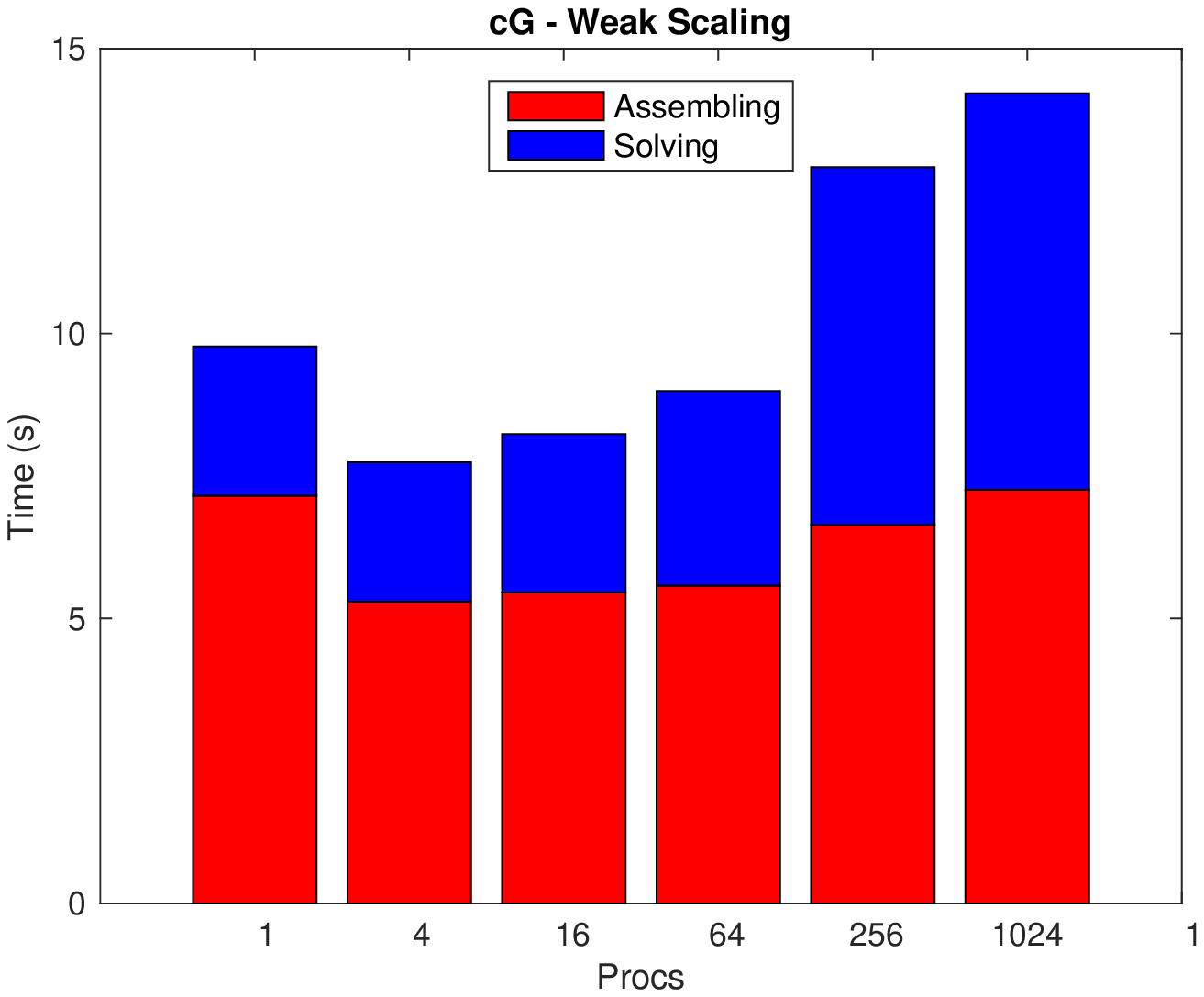}
  	\caption{$p=3$}
  \end{subfigure}
    \begin{subfigure}{0.328\textwidth}
        \includegraphics[width=0.99\textwidth]{\PathToPic 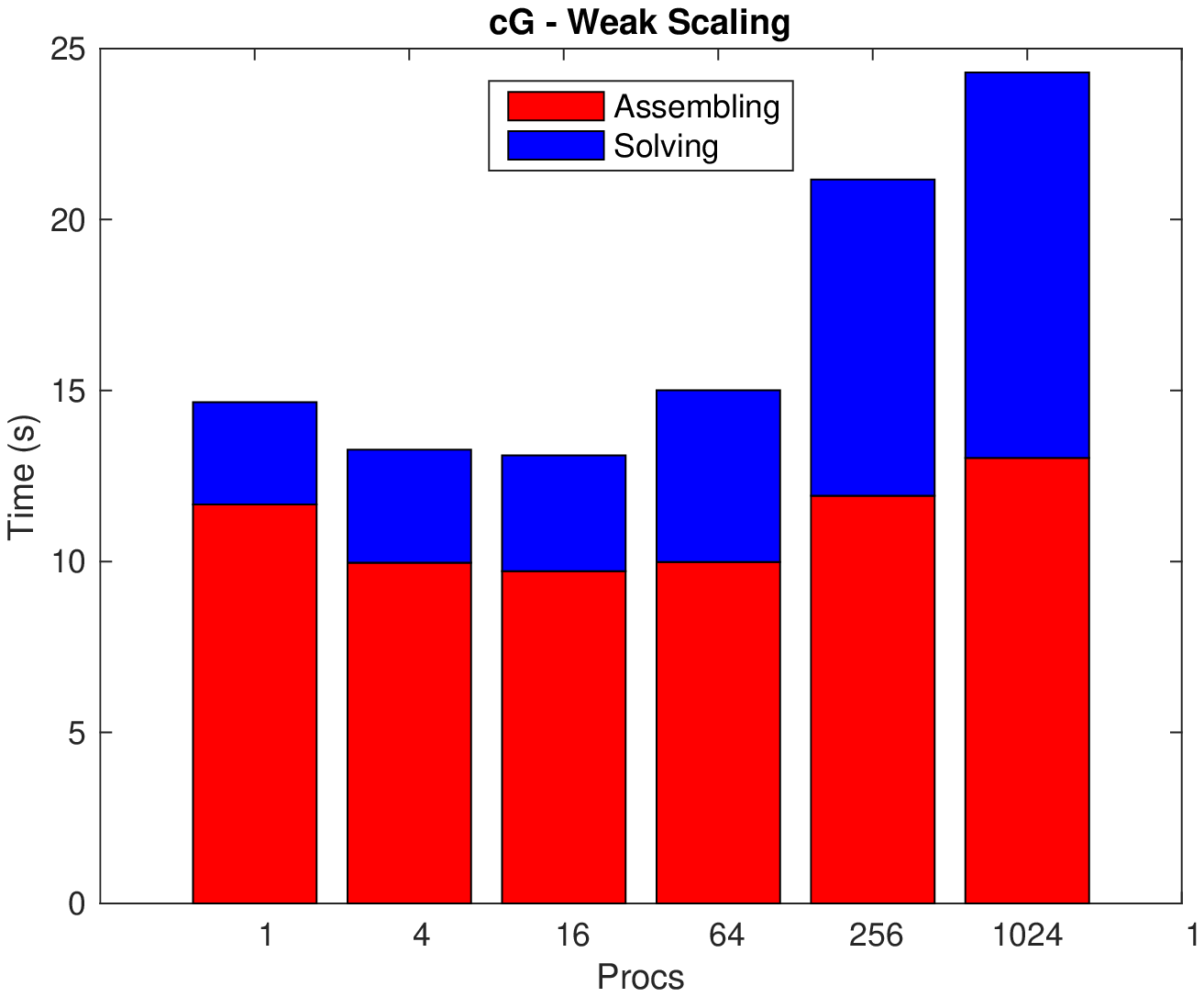}
  	\caption{$p=4$}
  \end{subfigure}
  \\
    \begin{subfigure}{0.328\textwidth}
        \includegraphics[width=0.99\textwidth]{\PathToPic 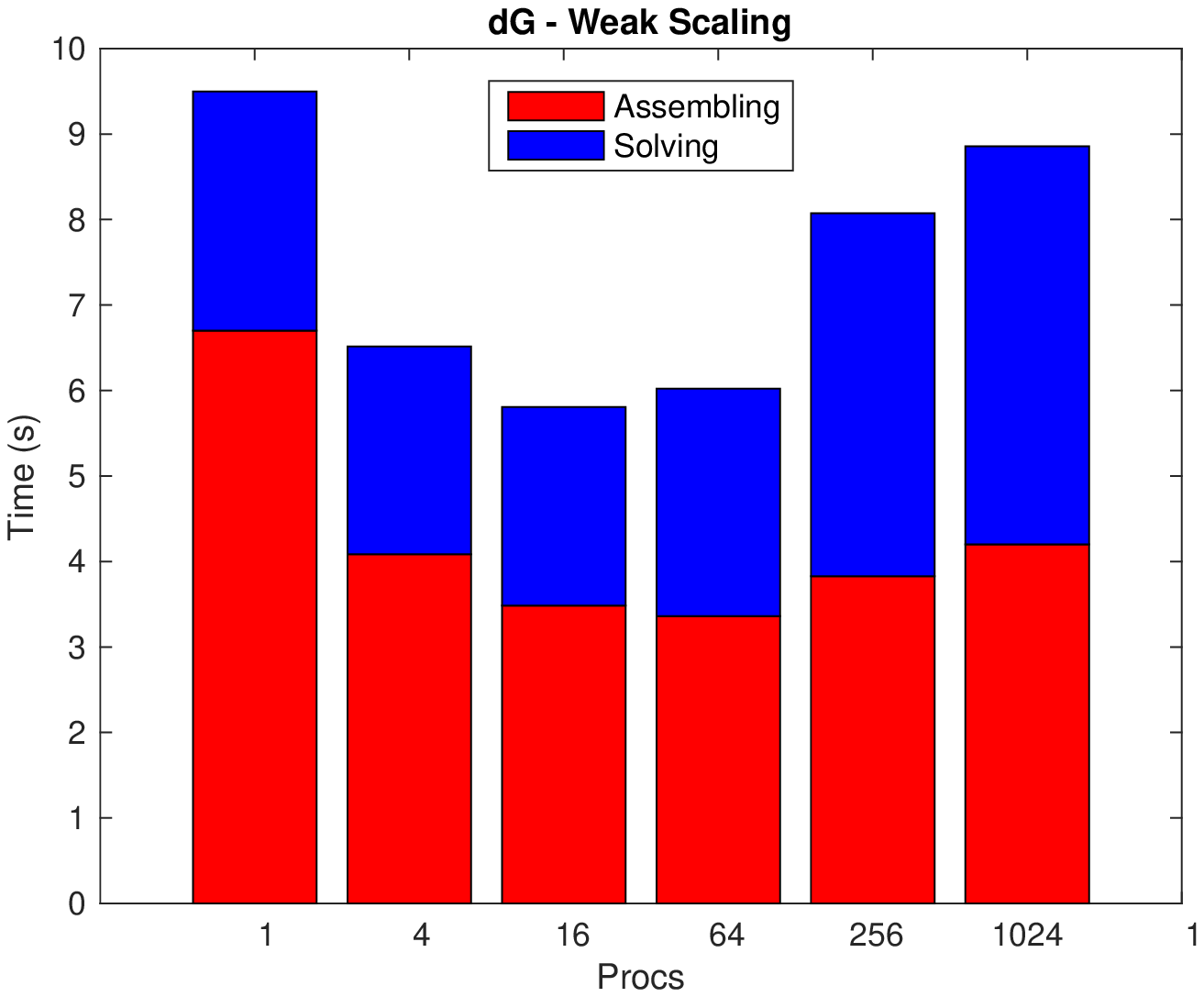}
  	\caption{$p=2$}
  \end{subfigure}
    \begin{subfigure}{0.328\textwidth}
        \includegraphics[width=0.99\textwidth]{\PathToPic 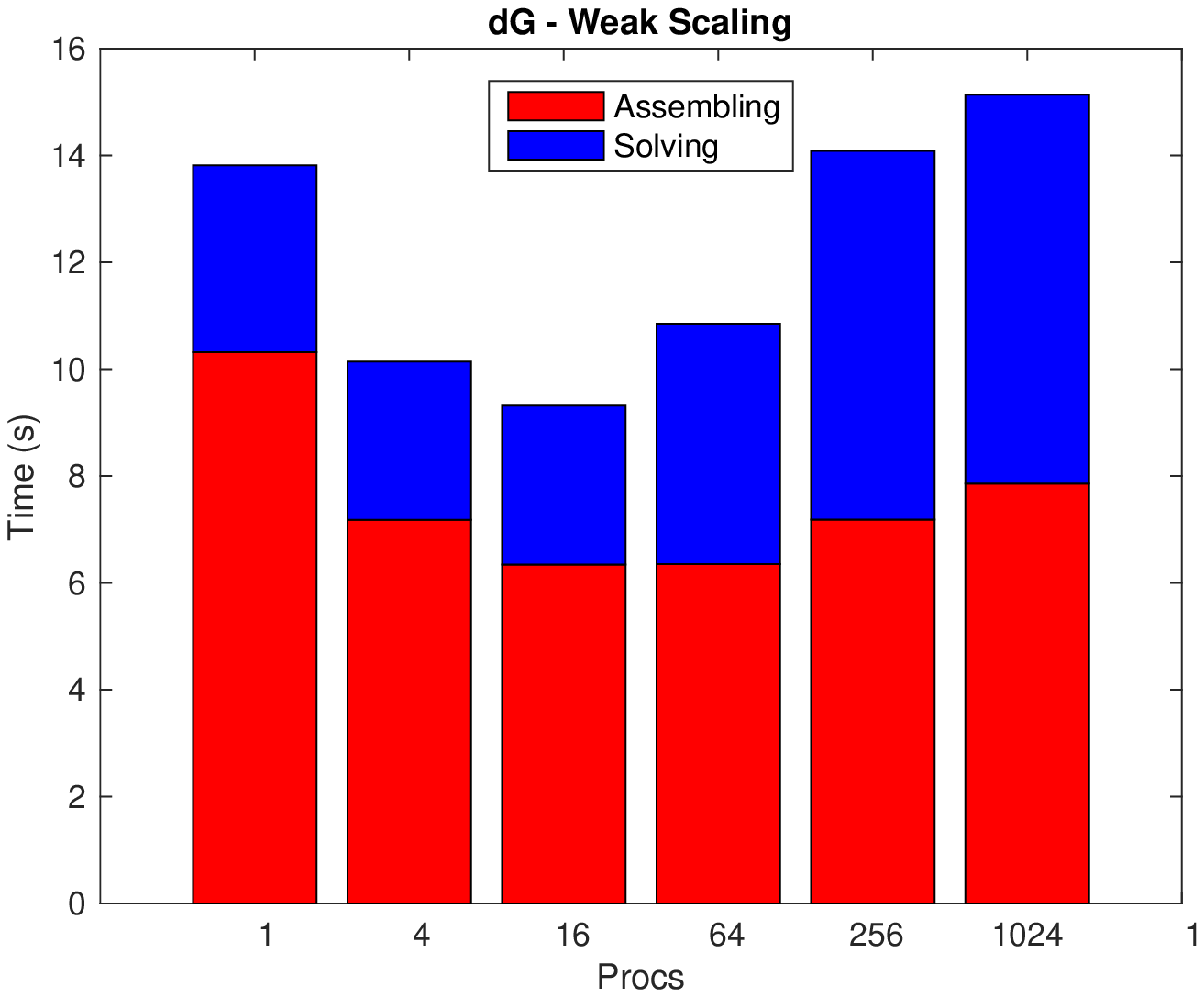}
  	\caption{$p=3$}
  \end{subfigure}
    \begin{subfigure}{0.328\textwidth}
        \includegraphics[width=0.99\textwidth]{\PathToPic 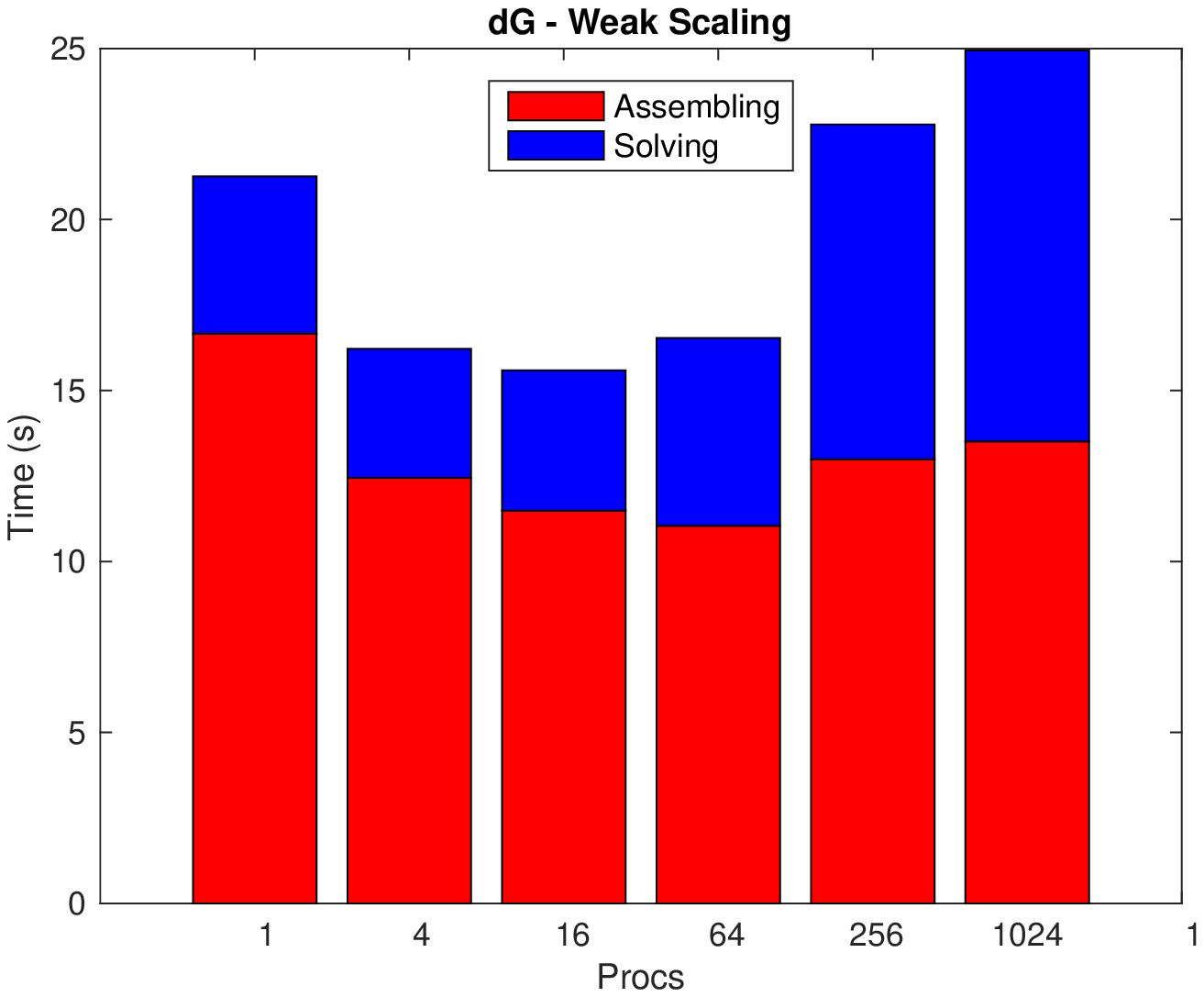}
  	\caption{$p=4$}
  \end{subfigure}
  \caption{Weak scaling of the cG-IETI-DP (first row) and dG-IETI-DP (second row) method for B-Spline degrees $p\in\{2,3,4\}$ in two dimensions. Each degree corresponds to one column. }
  \label{fig:Res_weak_2d}
\end{figure}

First, we consider the two dimensional case. We apply three initial refinements and start with a single processor and perform up to additional 5 refinements with maximum 1024 processors. We choose as primal variables continuous vertex values and edge averages. The results for degree $p\in\{2,3,4\}$ are illustrated in Figure~\ref{fig:Res_weak_2d}. The first row of figures corresponds to the cG-IETI-DP method, and the second one corresponds to the dG-IETI-DP method. The left column of Table~\ref{tab:weak_scaling1} summarizes timings and the speedup for the cG-IETI-DP method, whereas the right column presents the results for the dG-IETI-DP method.  For each method, we investigate the weak scaling for the assembling and solution phase. As in Figure~\ref{fig:Res_weak_2d}, we present the scaling and timings for $p\in\{2,3,4\}$.

\begin{table}  [h]
  \begin{footnotesize}
 \begin{tabular}{|r|r|c|c|c|c||r|c|c|c|c|}\hline
  & \multicolumn{2}{c|}{cG-IETI-DP}&\multicolumn{3}{c||}{ $p=2$} & \multicolumn{2}{c|}{dG-IETI-DP} &\multicolumn{3}{c|}{ $p=2$}        \\\hline
$\#$ procs &  $\#$dofs & Iter. & \makecell{Ass. \\ Time} & \makecell{Solv. \\ Time} &  \makecell{Total \\ Time}   &  $\#$dofs & Iter. & \makecell{Ass. \\ Time} & \makecell{Solv. \\ Time} &  \makecell{Total \\ Time}                 \\\hline
1          &        99104   &   7    &   4.8      &    1.9   &   6.7     &      133824   &  8    &    6.7  &   2.8   &      9.5                    \\\hline
4          &       328224   &   8   &    3.0     &     1.8  &    4.8    &       394688  &   9   &     4.1 &    2.4  &       6.5                   \\\hline
16         &      1179680   &   9    &   2.9      &    1.9   &   4.8     &     1309632   &  10   &    3.5  &   2.3   &      5.8                    \\\hline
64         &      4455456   &   10   &   3.0      &    2.4   &   5.4     &     4712384   &  11  &     3.4 &    2.7  &       6.1                   \\\hline 
256        &     17298464   &   11   &   3.5      &    4.2   &   7.7     &    17809344   &  11  &     3.8 &    4.2  &       8.0                   \\\hline 
1024       &     68150304   &   11  &    3.8     &     4.4  &    8.2    &     69169088  &   12  &     4.2 &    4.7  &       8.9                   \\\hline \hline
 & \multicolumn{2}{c|}{cG-IETI-DP}&\multicolumn{3}{c||}{ $p=3$}  &\multicolumn{2}{c|}{dG-IETI-DP}& \multicolumn{3}{c|}{ $p=3$}       \\\hline
$\#$ procs &  $\#$dofs & Iter. & \makecell{Ass. \\ Time} & \makecell{Solv. \\ Time} &  \makecell{Total \\ Time}   &  $\#$dofs & Iter. & \makecell{Ass. \\ Time} & \makecell{Solv. \\ Time} &  \makecell{Total \\ Time}                   \\\hline
1          &     120576   &   8    &    7.2   &  2.6  &      9.8      &      159264    &  8      &  10.3   &   3.5   &  13.8                    \\\hline
4          &     366080   &   9   &     5.3  &   2.4 &       7.7     &       436512   &   9      &   7.2   &   3.0   &  10.2                    \\\hline
16         &    1250304   &   10   &    5.5   &  2.8  &      8.3      &     1384224    &  10     &   6.3   &   3.0   &   9.3                    \\\hline
64         &    4591616   &   10   &    5.6   &  3.4  &      9.0      &     4852512    &  11     &   6.4   &   4.5   &  10.9                    \\\hline
256        &   17565696   &   11   &    6.6   &  6.3  &     12.9      &    18080544    &  12     &   7.2   &   6.9   &  14.1                    \\\hline 
1024       &   68679680   &   12  &     7.3  &   7.0 &      14.3     &     69702432   &   12    &    7.9  &    7.3  &   15.2                   \\\hline \hline
 & \multicolumn{2}{c|}{cG-IETI-DP}&\multicolumn{3}{c||}{ $p=4$}  & \multicolumn{2}{c|}{dG-IETI-DP}& \multicolumn{3}{c|}{ $p=4$}          \\\hline
$\#$ procs &  $\#$dofs & Iter. & \makecell{Ass. \\ Time} & \makecell{Solv. \\ Time} &  \makecell{Total \\ Time}   &  $\#$dofs & Iter. & \makecell{Ass. \\ Time} & \makecell{Solv. \\ Time} &  \makecell{Total \\ Time}                        \\\hline
1          &      144096   &   8    &    11.7  &    3.0&       14.7    &        186752 &  9       &  16.7  &     4.6  &    21.3                     \\\hline
4          &      405984   &   9    &    10.0  &    3.3 &      13.3     &       480384  &  10     &  12.4  &     3.8  &    16.2                     \\\hline
16         &     1322976   &   10   &     9.7    &  3.4   &    13.1       &    1460864  &  11     &  11.5  &     4.1  &    15.6                     \\\hline
64         &     4729824   &   11   &    10.0    &  5.0   &    15.0       &    4994688   &  11    &  11.0  &     5.5  &    16.5                     \\\hline
256        &    17834976   &   12  &     11.9   &   9.3  &     21.2      &    18353792   &   12   &  13.0  &     9.8  &    22.8                     \\\hline 
1024       &    69211104   &   13  &     13.0   &  11.3  &     24.3      &    70237824   &  13    &  13.5  &    11.4  &    24.9                     \\\hline
 \end{tabular}
 \end{footnotesize}
 \caption{Weak scaling results for the two dimensional testcase for the cG and dG IETI-DP method. Left column contains results for the cG variant and the right column for the dG version. Each row corresponds to a fixed B-Spline degree $p\in\{2,3,4\}$}
 \label{tab:weak_scaling1}
\end{table}

We observe that the time used for the assembling stays almost constant, hence shows quite optimal behaviour. However, the time for solving the system increases when refining and increasing the number of used processors. Especially, when considering the largest number of processors, we see a clear increase of the solution time. One reason is that the number of iterations slightly increases when increasing the system size. This is due to the quasi optimal condition number bound of the IETI-DP type methods, cf. \eqref{equ:kappa}. Secondly, as already pointed out in Section~\ref{sec:para}, the solving phase consists of more communication between processors, which cannot be completely overlapped with computations. Moreover, one also has to take in account global synchronization points in the conjugate gradient method.

\begin{figure}[h]
  \begin{subfigure}{0.328\textwidth}
        \includegraphics[width=0.99\textwidth]{\PathToPic 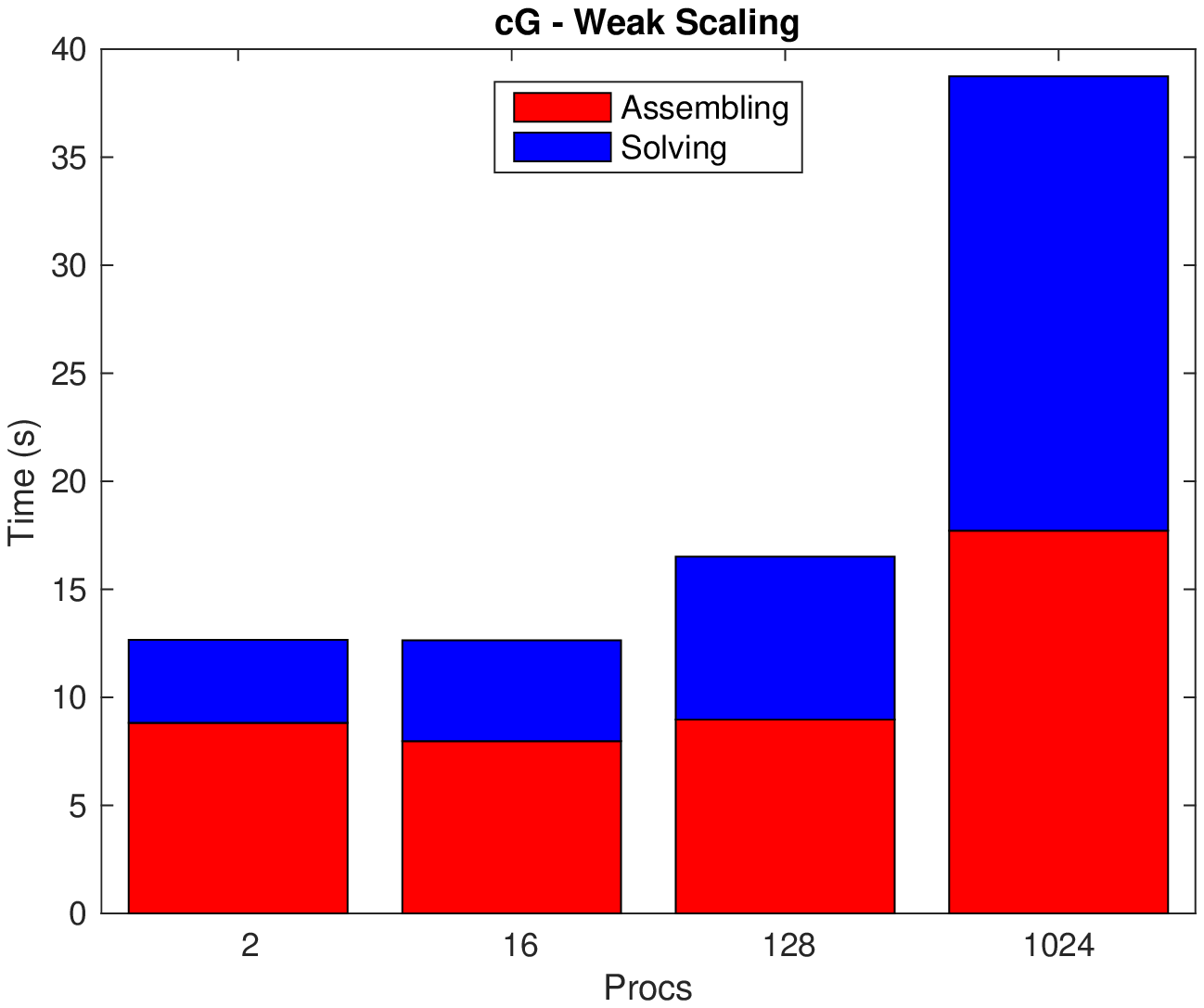}
  	\caption{$p=2$}
  \end{subfigure}
    \begin{subfigure}{0.328\textwidth}
        \includegraphics[width=0.99\textwidth]{\PathToPic 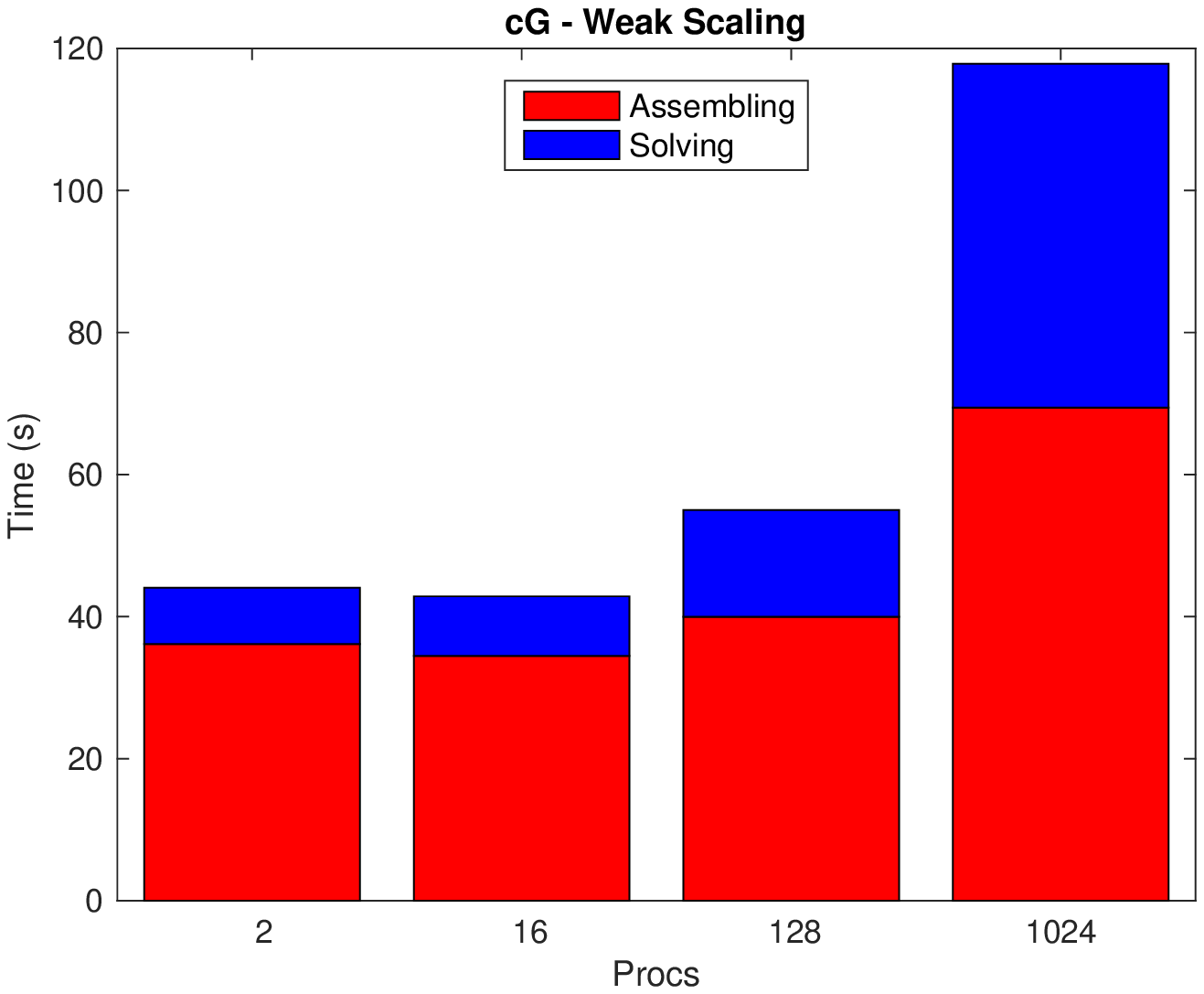}
  	\caption{$p=3$}
  \end{subfigure}
    \begin{subfigure}{0.328\textwidth}
        \includegraphics[width=0.99\textwidth]{\PathToPic 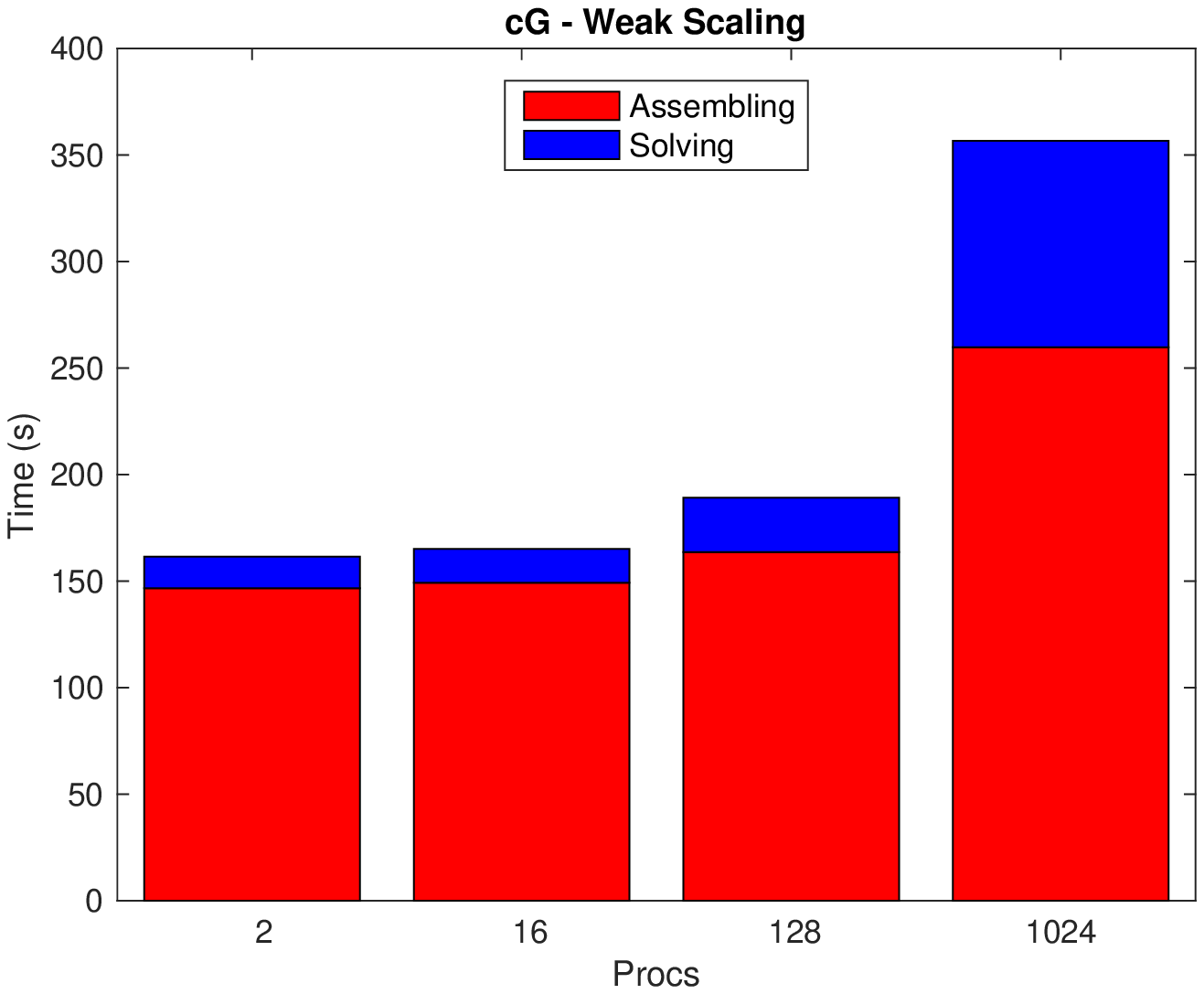}
  	\caption{$p=4$}
  \end{subfigure}
  \\
    \begin{subfigure}{0.328\textwidth}
        \includegraphics[width=0.99\textwidth]{\PathToPic 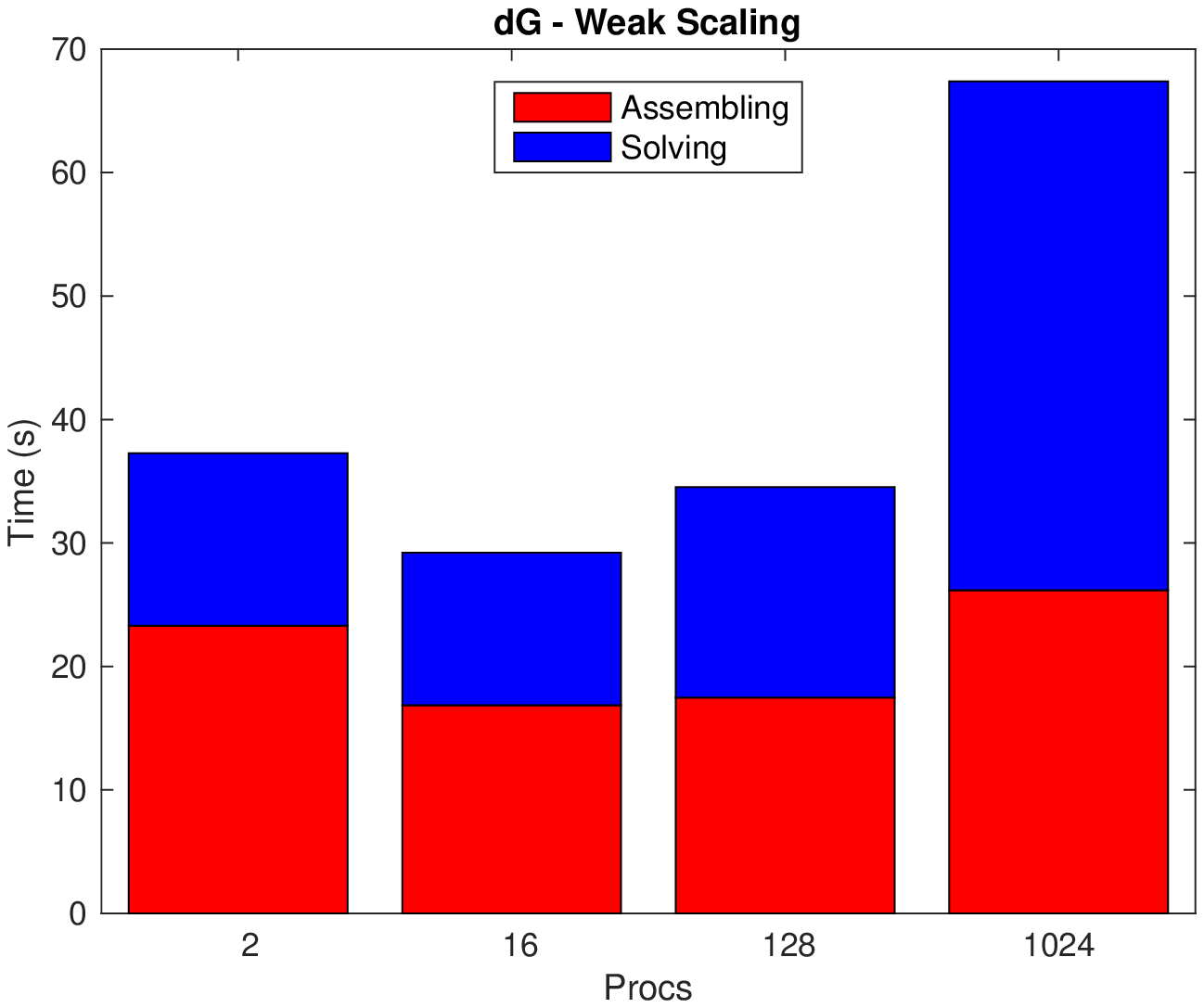}
  	\caption{$p=2$}
  \end{subfigure}
    \begin{subfigure}{0.328\textwidth}
        \includegraphics[width=0.99\textwidth]{\PathToPic 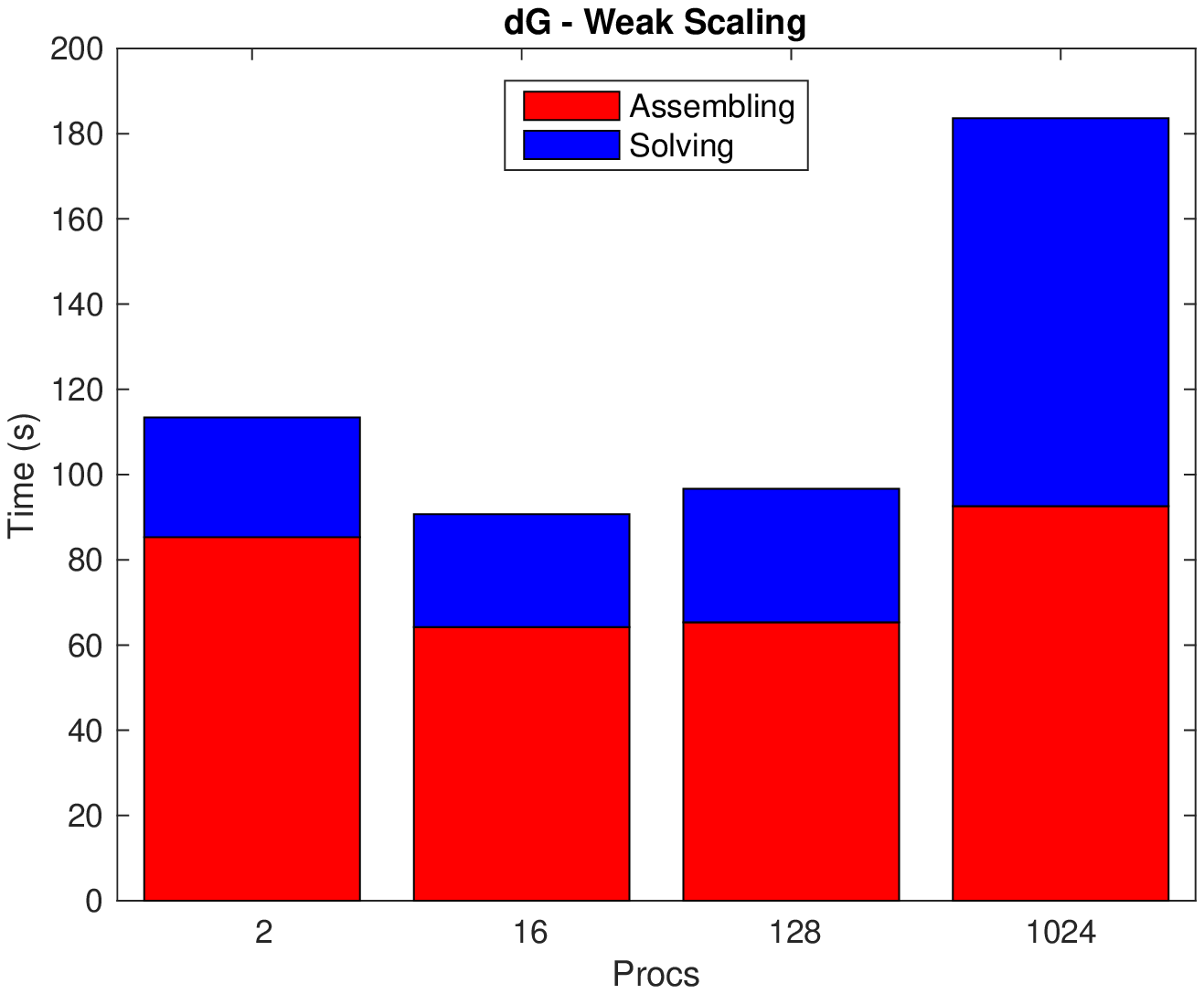}
  	\caption{$p=3$}
  \end{subfigure}
    \begin{subfigure}{0.328\textwidth}
        \includegraphics[width=0.99\textwidth]{\PathToPic 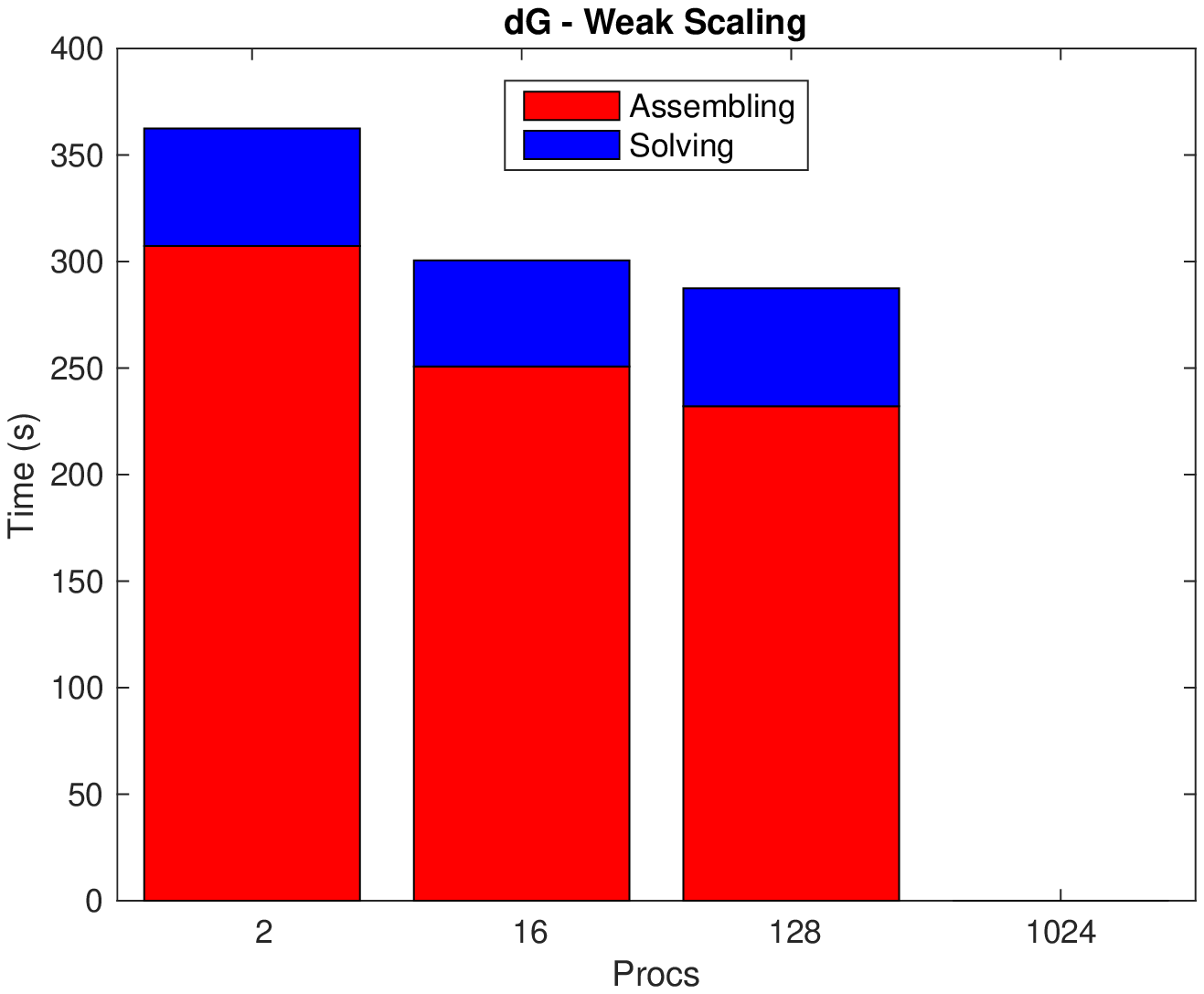}
  	\caption{$p=4$}
  \end{subfigure}
  \caption{Weak scaling of the cG-IETI-DP (first row) and dG-IETI-DP (second row) method for B-Spline degrees $p\in\{2,3,4\}$ in three dimensions. Each degree corresponds to one column. No timings are obtained in the case of 1024 cores in (f) due to memory limitations }
  \label{fig:Res_weak_3d}
\end{figure}

Next, we consider the weak scaling for the three dimensional case. As already indicated in the introduction of this section, we choose only continuous edge averages as primal variables. We perform the tests in the same way as for the two dimensional case. However, we already start with two processors and perform two initial refinements. Multiplying the number of used processors by 8 with each refinement, we end up again with 1024 processors on the finest grid. The two algorithms behave similar to the two dimensional case, where the assembling phase gives quite good results and the solver phase shows again an increasing time after each refinement. The results are visualized in Figure~\ref{fig:Res_weak_3d} and summarized in Table~\ref{tab:weak_scaling2}. Note, for the dG-IETI-DP method with $p=4$ and $\sim 54$~Mio. dofs, we exceeded the memory capacity of the cluster.

\begin{table} [h]
     \begin{footnotesize}
 \begin{tabular}{|r|r|c|c|c|c||r|c|c|c|c|}\hline
  & \multicolumn{2}{c|}{cG-IETI-DP}&\multicolumn{3}{c||}{ $p=2$}    & \multicolumn{2}{c|}{dG-IETI-DP} &\multicolumn{3}{c|}{ $p=2$}   \\\hline
$\#$ procs &  $\#$dofs & Iter. & \makecell{Ass. \\ Time} & \makecell{Solv. \\ Time} &  \makecell{Total \\ Time}   &  $\#$dofs & Iter. & \makecell{Ass. \\ Time} & \makecell{Solv. \\ Time} &  \makecell{Total \\ Time}                 \\\hline
2          &        220896   &   16    &    8.8  &     3.8    &     12.6     &        396932  & 26      &  23.3 &  14.0  &  37.3                    \\\hline                
16         &       1023200   &   17    &    8.0  &     4.7    &     12.7     &      1551400   & 27      &  16.9 &  12.3  &  29.2                    \\\hline                
128        &       5969376   &   17    &    9.0  &     7.5    &     16.5     &      7730288   & 28      &  17.5 &  17.0  &  34.5                    \\\hline                
1024       &      40238048   &   19    &   17.7  &    21.0    &     38.7     &     46577920   & 28      &  26.2 &  41.2  &  67.4                    \\\hline \hline         
 & \multicolumn{2}{c|}{cG-IETI-DP}&\multicolumn{3}{c||}{ $p=3$} &\multicolumn{2}{c|}{dG-IETI-DP}& \multicolumn{3}{c|}{ $p=3$}   \\\hline
$\#$ procs &  $\#$dofs & Iter. & \makecell{Ass. \\ Time} & \makecell{Solv. \\ Time} &  \makecell{Total \\ Time}   &  $\#$dofs & Iter. & \makecell{Ass. \\ Time} & \makecell{Solv. \\ Time} &  \makecell{Total \\ Time}                   \\\hline                
2          &       350840   &   17    &  36.2    &   7.9     &  44.1   &       598405  &  29     &  85.3  &   28.1    &  113.4                  \\\hline               
16         &      1361976   &   18   &   34.5   &    8.4    &   42.9   &      2005737  &  30     &  64.2  &   26.5    &   90.7                  \\\hline               
128        &      7020728   &   18   &   40.0   &   15.0    &   55.0   &      8985265  &  30     &  65.4  &   31.3    &   96.7                  \\\hline               
1024       &     43894200   &   21    &  69.4    &  48.4     & 117.8   &    50613825   &  31     &  92.6  &   91.0    &  183.6                   \\\hline \hline       
 & \multicolumn{2}{c|}{cG-IETI-DP}&\multicolumn{3}{c||}{ $p=4$}   &\multicolumn{2}{c|}{dG-IETI-DP}& \multicolumn{3}{c|}{ $p=4$}     \\\hline
$\#$ procs &  $\#$dofs     & Iter.   & \makecell{Ass. \\ Time}   & \makecell{Solv. \\ Time} &  \makecell{Total \\ Time} &  $\#$dofs    & Iter.   & \makecell{Ass. \\ Time}& \makecell{Solv. \\ Time}&  \makecell{Total \\ Time}                        \\\hline           
2          &      523776   &  18     &    146.6    &   14.8     &   161.4  &        853878   &  32     &  307.3   &   55.1    & 362.4                    \\\hline
16         &     1768320   &  18     &    149.2    &   15.9     &   165.1  &       2538650   &  32     &  250.7   &   49.8    & 300.5                    \\\hline
128        &     8188800   &  20     &    163.6    &   25.6     &   189.2  &      10367970   &  34     &  232.0   &   55.4    & 287.4                    \\\hline
1024       &    47765376   &  22     &    259.7    &   96.9     &   356.6  &  $\sim$54000000 &  x      &  x       &     x     &  x                     \\\hline
 \end{tabular}
    \end{footnotesize}
 \caption{Weak scaling results for the three dimensional testcase for the cG and dG IETI-DP method. Left column contains results for the cG variant and the right column for the dG version. Each row corresponds to a fixed B-Spline degree $p\in\{2,3,4\}$. No timings are available for the dG-IETI-DP method with $p=4$ on 1024 cores due to memory limitations. }
 \label{tab:weak_scaling2}
\end{table}

\subsection{Strong scaling}
\label{sec:strong}
Secondly, we are investigating the strong scaling behaviour. Now we fix the problem size and increase the number of processors. In the optimal case, the time used by a certain quantity reduces in the same way as the number of used processors increases. We use the same primal variables for the strong scaling studies as in the weak scaling studies in Section~\ref{sec:weak}.

\begin{figure}[h]
  \begin{subfigure}{0.48\textwidth}
        \includegraphics[width=0.99\textwidth]{\PathToPic 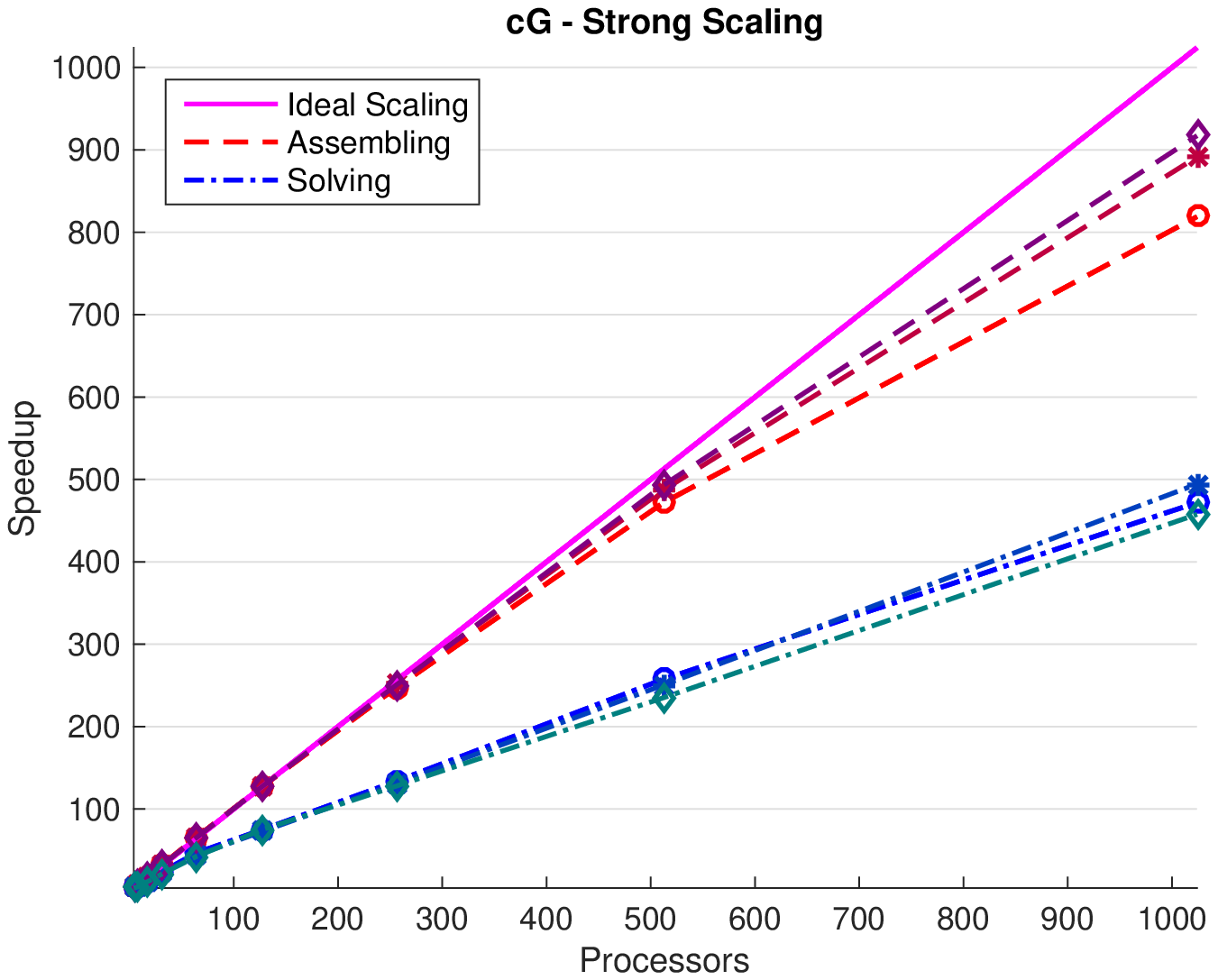}
  	\caption{cG-IETI-DP}
  \end{subfigure}
    \begin{subfigure}{0.48\textwidth}
        \includegraphics[width=0.99\textwidth]{\PathToPic 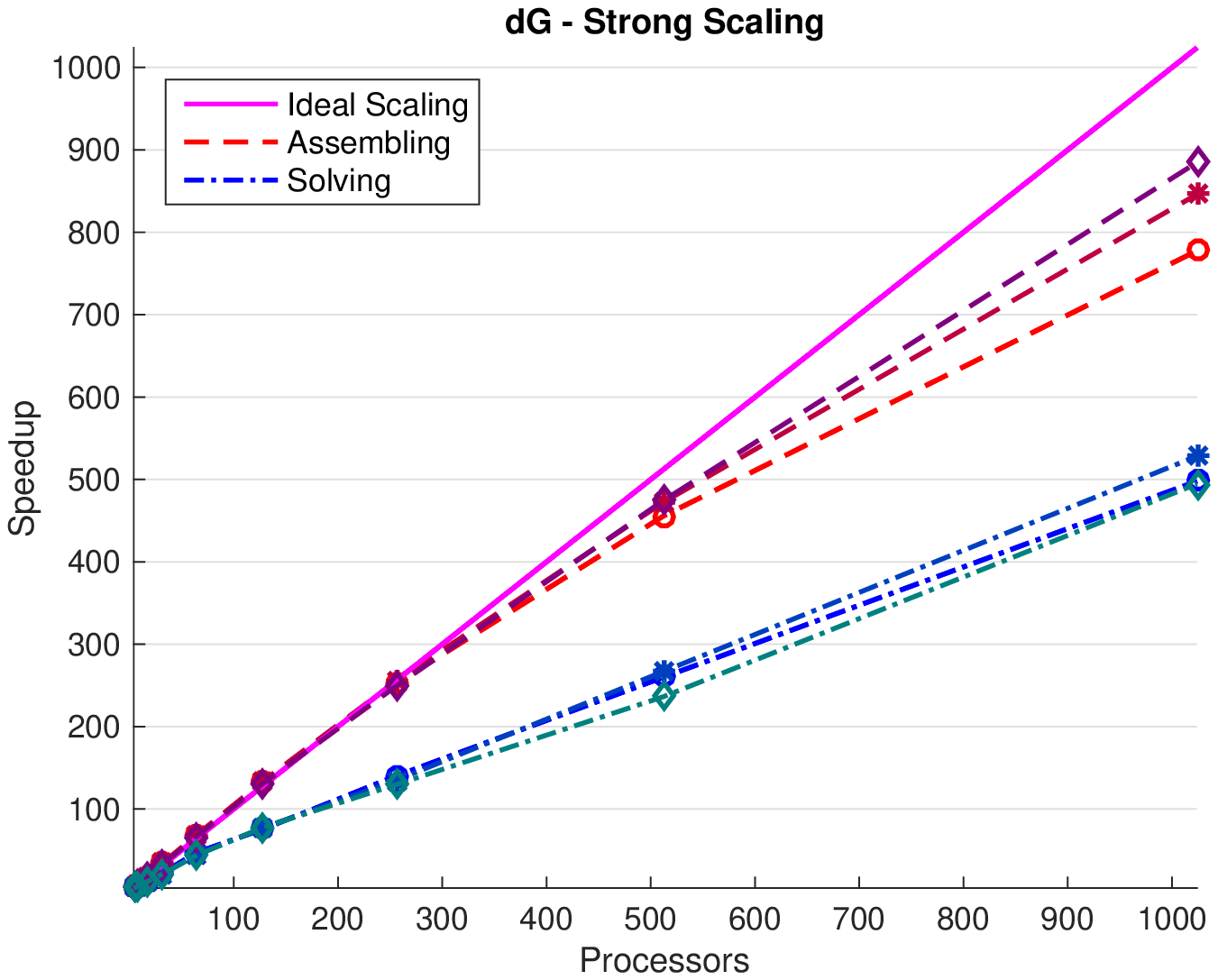}
  	\caption{dG-IETI-DP}
  \end{subfigure}
  \caption{Strong scaling of the cG-IETI-DP (left column) and dG-IETI-DP (right column) method for B-Spline degrees $p\in\{2,3,4\}$ in two dimensions.  The markers $\{\circ,*,\Diamond\}$ as well as different shades of red (assembling phase) and blue (solver phase) correspond to the degrees $\{2,3,4\}$.}
  \label{fig:strong_2d}
\end{figure}

Again as in Section~\ref{sec:weak}, we begin with the two dimensional example. We perform $7$ initials refinements and end up with $17$~Mio. dofs on 1024 subdomains. We start already with $4$ processors in the initial case and do $8$ refinements until we reach 1024 cores. Similar to Section~\ref{sec:weak}, the results for $p\in\{2,3,4\}$ are illustrated in Figure~\ref{fig:strong_2d} and summarized in Table~\ref{tab:strong_scaling1}.             

\begin{table}   [h]
     \begin{footnotesize}
 \begin{tabular}{|r|c|c|c|c||c|c|c|c||c|c|c|c|}\hline
2d  & \multicolumn{4}{c||}{ $p=2$} & \multicolumn{4}{c||}{$p=3$}& \multicolumn{4}{c|}{$p=4$}\\ \hline
   \tiny{cG-IETI-DP}  &  \multicolumn{2}{c|}{\makecell{assembling \\ phase}} & \multicolumn{2}{c||}{\makecell{solving \\ phase}} &  \multicolumn{2}{c|}{\makecell{assembling \\ phase}} & \multicolumn{2}{c||}{\makecell{solving \\ phase}}&   \multicolumn{2}{c|}{\makecell{assembling \\ phase}}& \multicolumn{2}{c|}{\makecell{solving \\ phase}} \\ \hline
$\#$ procs & Time & Sp. & Time & Sp. & Time & Sp. & Time & Sp. & Time & Sp. & Time & Sp.                     \\ \hline
    4  &  190.8 &   4 & 138.2 &   4   &  364.5 &   4 & 202.6  &   4   &  653.7 &   4 & 279.1 &   4            \\ \hline
    8  &   94.8 &   8 &  88.8 &   6   &  181.2 &   8 & 141.8  &   6   &  325.9 &   8 & 202.7 &   6            \\ \hline
  16   &   47.1 &  16 &  45.9 &  12   &   89.9 &  16 &  71.9  &  11   &  162.3 &  16 & 102.6 &  11            \\ \hline                                  
  32   &   23.1 &  32 &  22.9 &  24   &   44.5 &  32 &  35.7  &  23   &   80.4 &  32 &  51.3 &  22            \\ \hline                                  
  64   &   11.6 &  65 &  11.8 &  46   &   22.4 &  65 &  18.5  &  44   &   40.2 &  64 &  26.3 &  42            \\ \hline                                  
  128  &    5.9 & 127 &   7.3 &  75   &   11.3 & 128 &  11.1  &  73   &   20.4 & 128 &  14.8 &  75            \\ \hline                                  
  256  &    3.0 & 247 &   4.1 & 133   &    5.7 & 251 &   6.1  & 131   &   10.4 & 250 &   8.7 & 128            \\ \hline                                  
  512  &    1.6 & 471 &   2.1 & 257   &    2.9 & 487 &   3.2  & 250   &    5.3 & 493 &   4.7 & 235            \\ \hline                                  
  1024 &    0.9 & 819 &   1.1 & 472   &    1.6 & 891 &   1.6  & 494   &    2.8 & 917 &   2.4 & 456            \\ \hline\hline                                                
  \tiny{dG-IETI-DP}  &  \multicolumn{2}{c|}{\makecell{assembling \\ phase}} & \multicolumn{2}{c||}{\makecell{solving \\ phase}} &  \multicolumn{2}{c|}{\makecell{assembling \\ phase}} & \multicolumn{2}{c||}{\makecell{solving \\ phase}}&  \multicolumn{2}{c|}{\makecell{assembling \\ phase}}& \multicolumn{2}{c|}{\makecell{solving \\ phase}} \\ \hline
  $\#$ procs  & Time & Sp. & Time & Sp.   & Time & Sp. & Time & Sp.     & Time & Sp. & Time & Sp.             \\ \hline
   4    & 216.6&    4& 144.0 &   4     & 402.2&   4& 225.5 &    4  & 711.8&    4& 294.2 &    4           \\ \hline
   8    & 106.9&    8&  92.9 &   6     & 199.5&   8& 156.8 &    5  & 352.6&    8& 210.6 &    5           \\ \hline
 16     &  52.5&   16&  47.9 &  12     &  98.1&  16&  80.4 &   11  & 174.3&   16& 106.4 &   11           \\ \hline
 32     &  25.1&   34&  25.7 &  22     &  47.5&  33&  42.2 &   21  &  84.9&   33&  55.2 &   21           \\ \hline
 64     &  12.7&   68&  12.2 &  47     &  23.9&  67&  20.4 &   44  &  42.9&   66&  27.0 &   43           \\ \hline
 128    &   6.5&  132&   7.6 &  75     &  12.0& 134&  11.6 &   77  &  21.7&  131&  15.1 &   77           \\ \hline
 256    &   3.4&  252&   4.1 & 140     &   6.2& 255&   6.6 &  135  &  11.4&  249&   9.0 &  129           \\ \hline
 512    &   1.9&  455&   2.2 & 260     &   3.4& 472&   3.3 &  267  &   6.0&  474&   4.9 &  236           \\ \hline
 1024   &   1.1&  777&   1.1 & 498     &   1.9& 846&   1.7 &  528  &   3.2&  885&   2.3 &  494           \\ \hline
 \end{tabular}
   \end{footnotesize}
    \caption{Strong scaling results: Time (s) and Speedup for $p\in\{2,3,4\}$ in two dimensions having approximately 17 Mio. dofs. First row shows results for the cG variant of the IETI-DP method, whereas the second row contains results for the dG version. Each column corresponds to a degree $p$.  }
    \label{tab:strong_scaling1}
\end{table}

We observe that the assembling phase has a quite good scaling performance, as already observed for the weak scaling results in Section~\ref{sec:weak}. Moreover, the higher the B-Spline degree, the better the parallel performance behaves. This holds due to increased computational costs for the parallel part. Similar to the weak scaling results, the solver phase does not provide such an excellent scaling as the assembling phase. Still, we obtain a scaling from around 500 when using 1024 processors. We note that the degree of the B-Splines does not seem to have such a significant effect on the scaling for the solver phase as for the assembling phase. 

\begin{figure}[h]
    \begin{subfigure}{0.48\textwidth}
        \includegraphics[width=0.99\textwidth]{\PathToPic 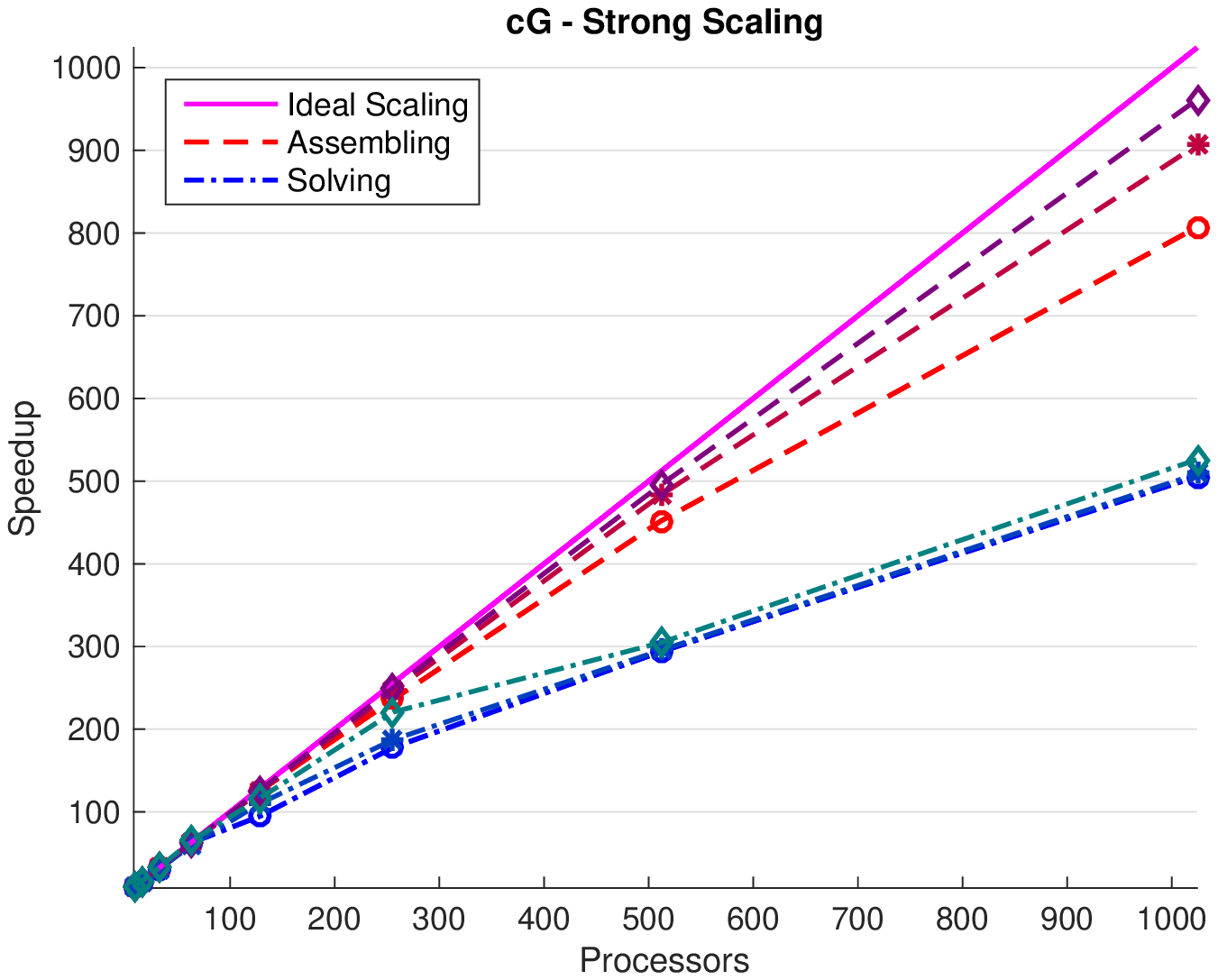}
  	\caption{cG-IETI-DP}
  \end{subfigure}
    \begin{subfigure}{0.48\textwidth}
        \includegraphics[width=0.99\textwidth]{\PathToPic 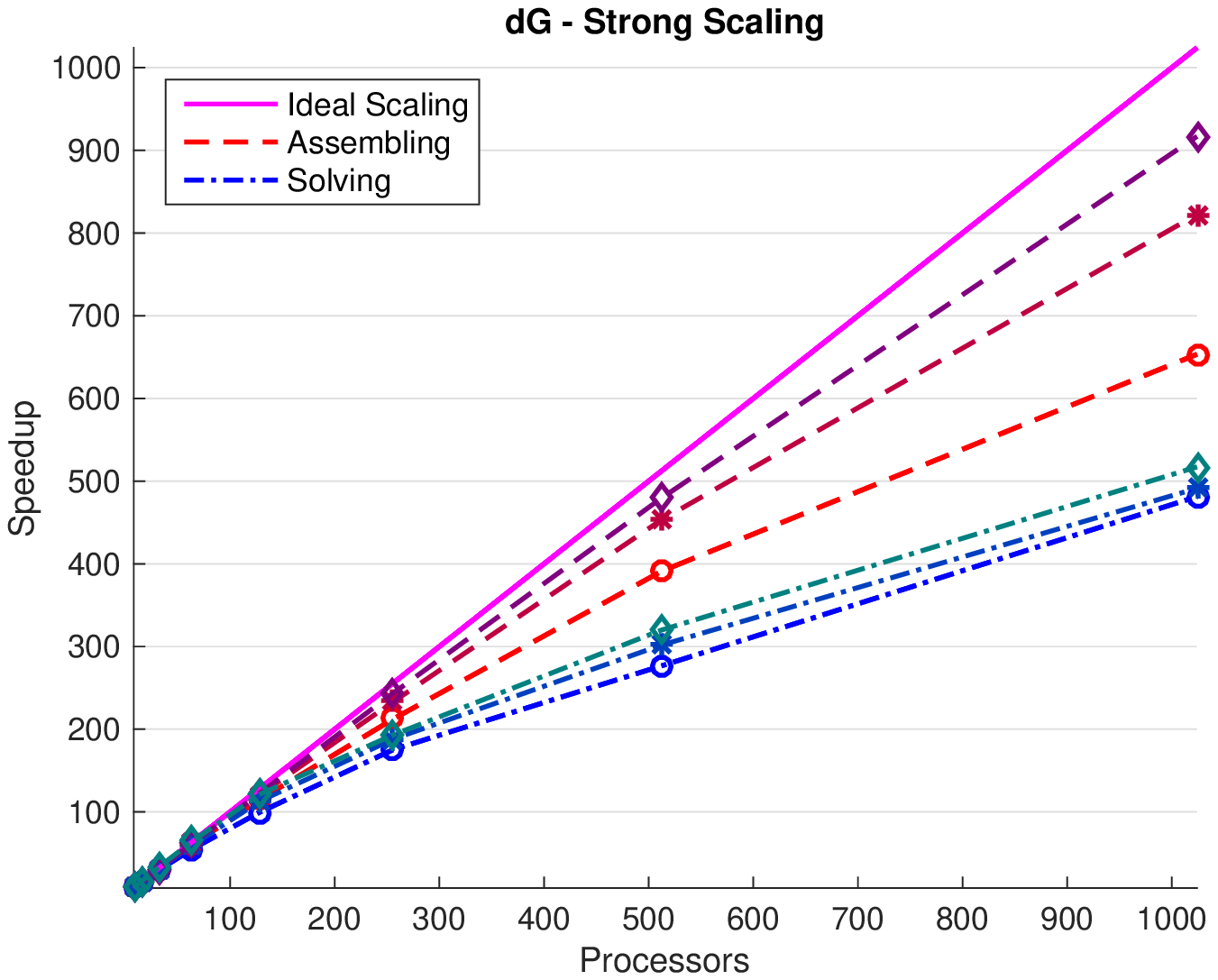}
  	\caption{dG-IETI-DP}
  \end{subfigure}
  \caption{Strong scaling of the cG-IETI-DP (left column) and dG-IETI-DP (right column) method for B-Spline degrees $p\in\{2,3,4\}$ in three dimensions.  The markers $\{\circ,*,\Diamond\}$ as well as different shades of red (assembling phase) and blue (solver phase) correspond to the degrees $\{2,3,4\}$.}
  \label{fig:strong_3d}
\end{figure}

In the three dimensional example we perform four initial refinements and obtain around 5 Mio. dofs. The presentation of the results is done in the same way as in the previous examples, see Figure~\ref{fig:strong_3d} and Table~\ref{tab:strong_scaling2}. Also in three dimensions the cG-IETI-DP algorithms behaves very similarly to the two dimensional case,   showing excellent scaling results. However, the dG version of the algorithm shows a good scaling but not as promising as cG version. Especially, when considering $p=2$, we observe degraded scalability for the assembling phase. Having a closer look at the timings, we observe that this originates from small load imbalances in the interior domains, due to the additional layer of dofs and the larger number of primal variables. The latter one leads to an increased time in solving \eqref{HL:equ:KC_basis}, due to a larger number of right hand sides on the interior subdomains.  One can further optimize the three dimensional case, by considering different strategies for the primal variables, where one aims for smaller and more equally distributed numbers of primal variables.

\begin{table}  [h]
    \begin{footnotesize}
 \begin{tabular}{|r|c|c|c|c||c|c|c|c||c|c|c|c|}\hline
  3d& \multicolumn{4}{c||}{ $p=2$} & \multicolumn{4}{c||}{$p=3$}& \multicolumn{4}{c|}{$p=4$}\\ \hline
   \tiny{cG-IETI-DP}  &  \multicolumn{2}{c|}{\makecell{assembling \\ phase}} & \multicolumn{2}{c||}{\makecell{solving \\ phase}} &  \multicolumn{2}{c|}{\makecell{assembling \\ phase}} & \multicolumn{2}{c||}{\makecell{solving \\ phase}}&   \multicolumn{2}{c|}{\makecell{assembling \\ phase}}& \multicolumn{2}{c|}{\makecell{solving \\ phase}} \\ \hline
$\#$ procs & Time & Sp. & Time & Sp. & Time & Sp. & Time & Sp. & Time & Sp. & Time & Sp.                   \\ \hline
    8  &     140.4   &     8   &  93.7   &    8  &  624.8  &     8 &   205.3  &     8  &  2565.7  &     8   &  382.6   &     8    \\ \hline
  16   &      70.5   &    16   &  47.7   &   16  &  312.8  &    16 &   103.4  &    16  &  1285.6  &    16   &  182.5   &    17    \\ \hline                                  
  32   &      35.2   &    32   &  25.0   &   30  &  157.0  &    32 &    53.7  &    31  &   643.3  &    32   &   93.7   &    33    \\ \hline                                  
  64   &      17.5   &    64   &  11.9   &   63  &   78.6  &    64 &    27.2  &    60  &   322.5  &    64   &   46.6   &    66    \\ \hline                                  
  128  &       9.0   &   125   &   7.9   &   95  &   40.2  &   124 &    14.9  &   110  &   163.0  &   126   &   26.4   &   116    \\ \hline                                  
  256  &       4.8   &   236   &   4.2   &  178  &   20.4  &   245 &     8.8  &   187  &    82.1  &   250   &   13.9   &   221    \\ \hline                                  
  512  &       2.5   &   452   &   2.6   &  294  &   10.3  &   483 &     5.6  &   296  &    41.4  &   496   &   10.0   &   305    \\ \hline                                  
  1024 &       1.4   &   807   &   1.5   &  506  &    5.5  &   906 &     3.2  &   509  &    21.4  &   961   &    5.8   &   526    \\ \hline\hline                                                
  \tiny{dG-IETI-DP} &  \multicolumn{2}{c|}{\makecell{assembling \\ phase}} & \multicolumn{2}{c||}{\makecell{solving \\ phase}} &  \multicolumn{2}{c|}{\makecell{assembling \\ phase}} & \multicolumn{2}{c||}{\makecell{solving \\ phase}}&  \multicolumn{2}{c|}{\makecell{assembling \\ phase}}& \multicolumn{2}{c|}{\makecell{solving \\ phase}} \\ \hline
  $\#$ procs  & Time & Sp. & Time & Sp.   & Time & Sp. & Time & Sp.     & Time & Sp. & Time & Sp.             \\ \hline      
 8      &   249.6   &     8   &  210.8    &    8   &    985.8    &    8   & 433.2    &     8   &   3588.1    &    8   &  854.8    &     8          \\ \hline  
 16     &   126.2   &    16   &  106.6    &   16   &    498.6    &   16   & 217.1    &    16   &   1792.9    &   16   &  405.7    &    17           \\ \hline                                              
 32     &    65.0   &    31   &   56.6    &   30   &    255.4    &   31   & 110.3    &    31   &    913.5    &   31   &  205.0    &    33           \\ \hline                                              
 64     &    33.1   &    60   &   30.5    &   55   &    128.5    &   61   &  58.2    &    60   &    460.0    &   62   &  105.9    &    65           \\ \hline                                              
 128    &    17.4   &   115   &   17.0    &   99   &     65.5    &  120   &  30.7    &   113   &    234.3    &  123   &   56.4    &   121           \\ \hline                                              
 256    &     9.4   &   212   &    9.6    &  175   &     33.8    &  233   &  18.4    &   188   &    117.8    &  244   &   35.4    &   193           \\ \hline                                              
 512    &     5.1   &   391   &    6.1    &  277   &     17.4    &  453   &  11.5    &   302   &     59.9    &  479   &   21.4    &   320            \\ \hline                                              
 1024   &     3.1   &   653   &    3.5    &  481   &      9.6    &  822   &   7.1    &   491   &     31.3    &  917   &   13.2    &   517            \\ \hline                                              
 \end{tabular}  
  \end{footnotesize}
    \caption{Strong scaling results: Time (s) and Speedup for $p\in\{2,3,4\}$ in three dimensions having approximately 5 Mio. dofs. First row shows results for the cG variant of the IETI-DP method, whereas the second row contains results for the dG version. Each column corresponds to a degree $p$.  }
    \label{tab:strong_scaling2}
\end{table}                                                                                                              

\subsection{Study on the number of $S_{\Pi\Pi}^{-1}$ holders}                                                                                              
\label{sec:diffHolder}
In this last section of the numerical experiments, we want to investigate the influence of the number of holders of $S_{\Pi\Pi}^{-1}$ on the scaling behaviour. As already indicated in Section~\ref{sec:para_AccDist}, if more processors hold the LU-factorization of the coarse grid matrix, it is possible to decrease the communication effort after applying $S_{\Pi\Pi}^{-1}$, while having more communication before the application. The advantage of this strategy is to be able to have a better overlap of communication with computations. However one has to take into account, that this also increases the communication in the assembling phase, since the local contribution $S\sMP_{\Pi\Pi}$ has to be sent to all the master processors.

\begin{figure}[h]
  \begin{subfigure}{0.328\textwidth}
        \includegraphics[width=0.99\textwidth]{\PathToPic 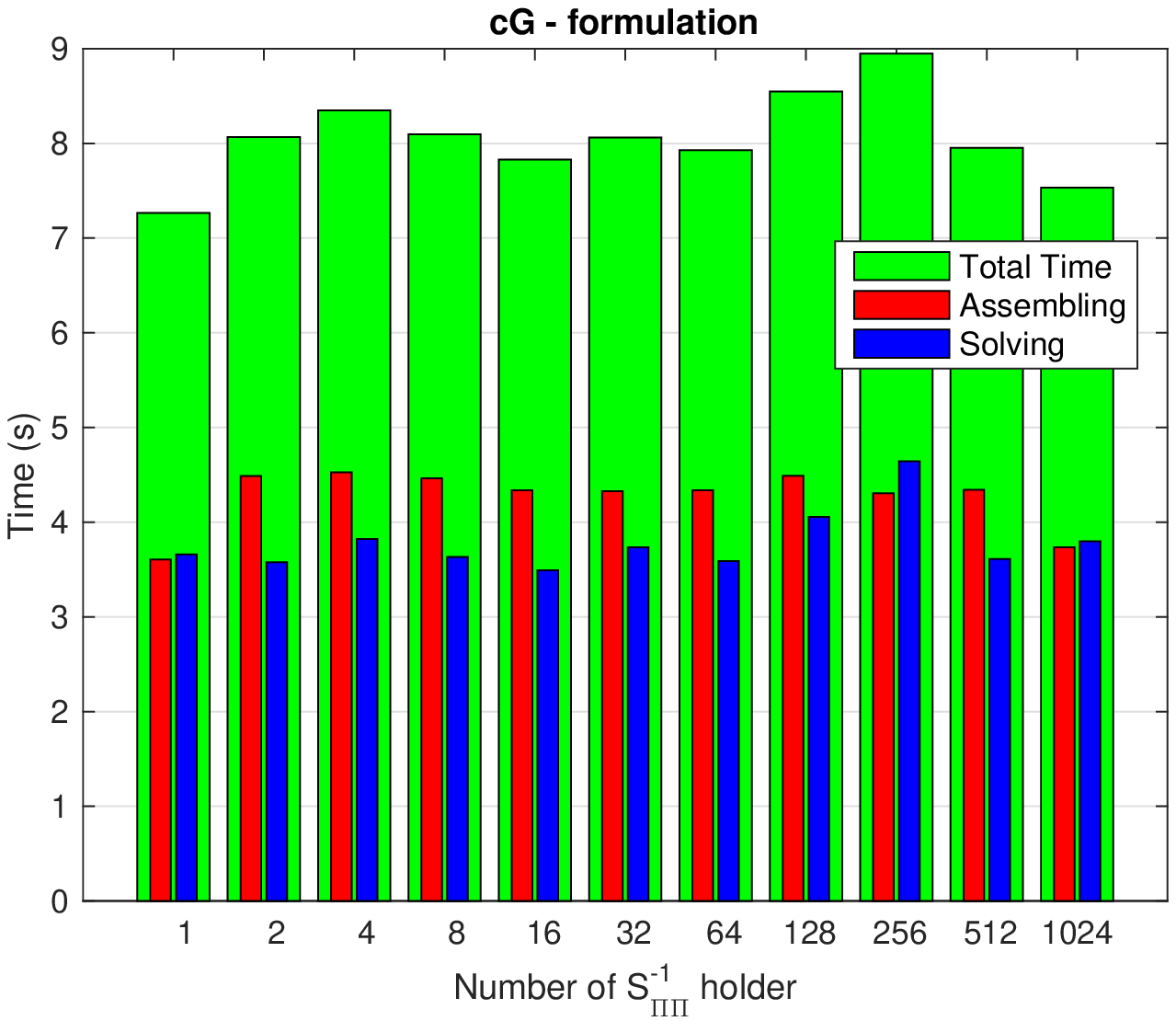}
  	\caption{$p=2$}
  \end{subfigure}
    \begin{subfigure}{0.328\textwidth}
        \includegraphics[width=0.99\textwidth]{\PathToPic 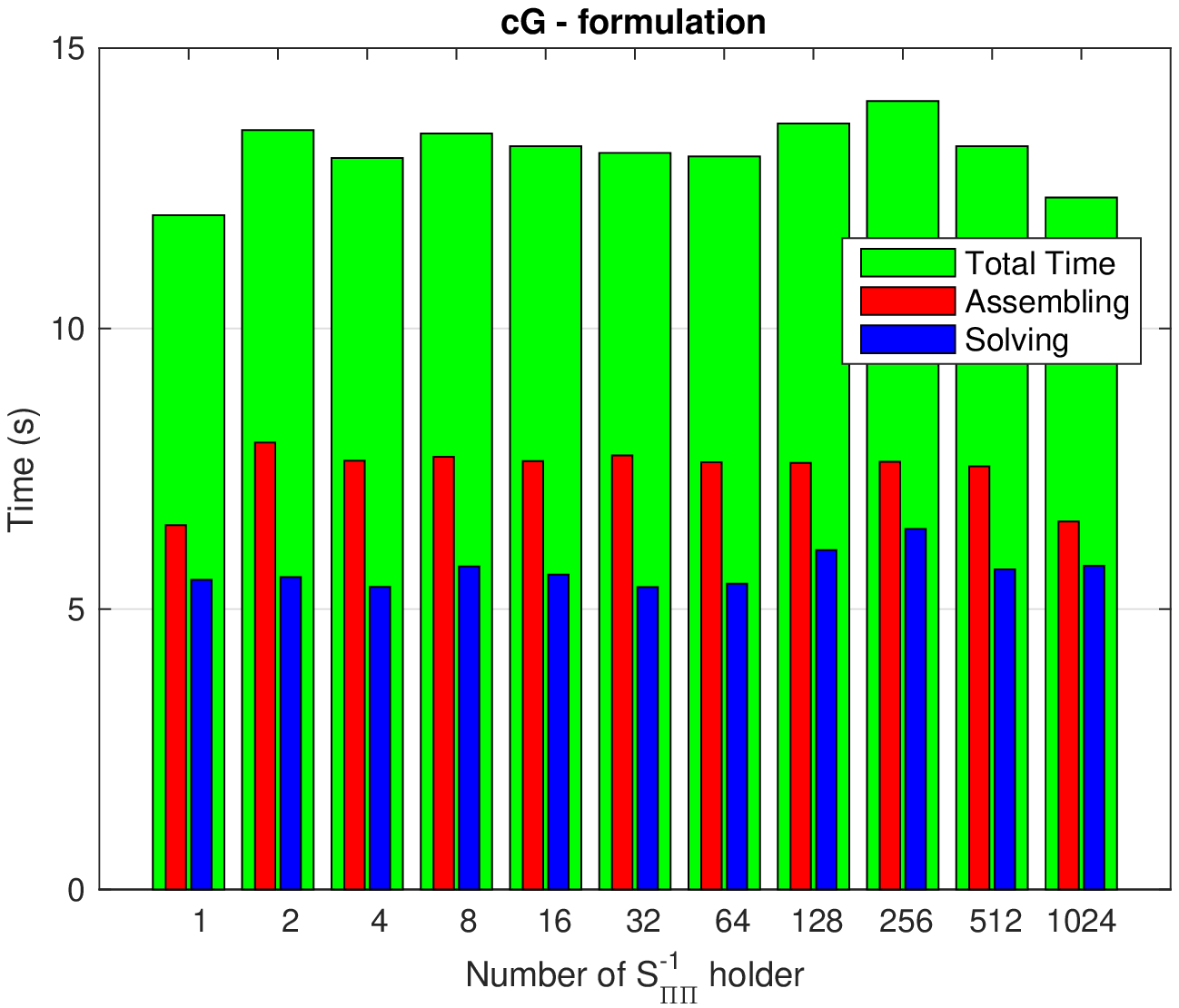}
  	\caption{$p=3$}
  \end{subfigure}
    \begin{subfigure}{0.328\textwidth}
        \includegraphics[width=0.99\textwidth]{\PathToPic 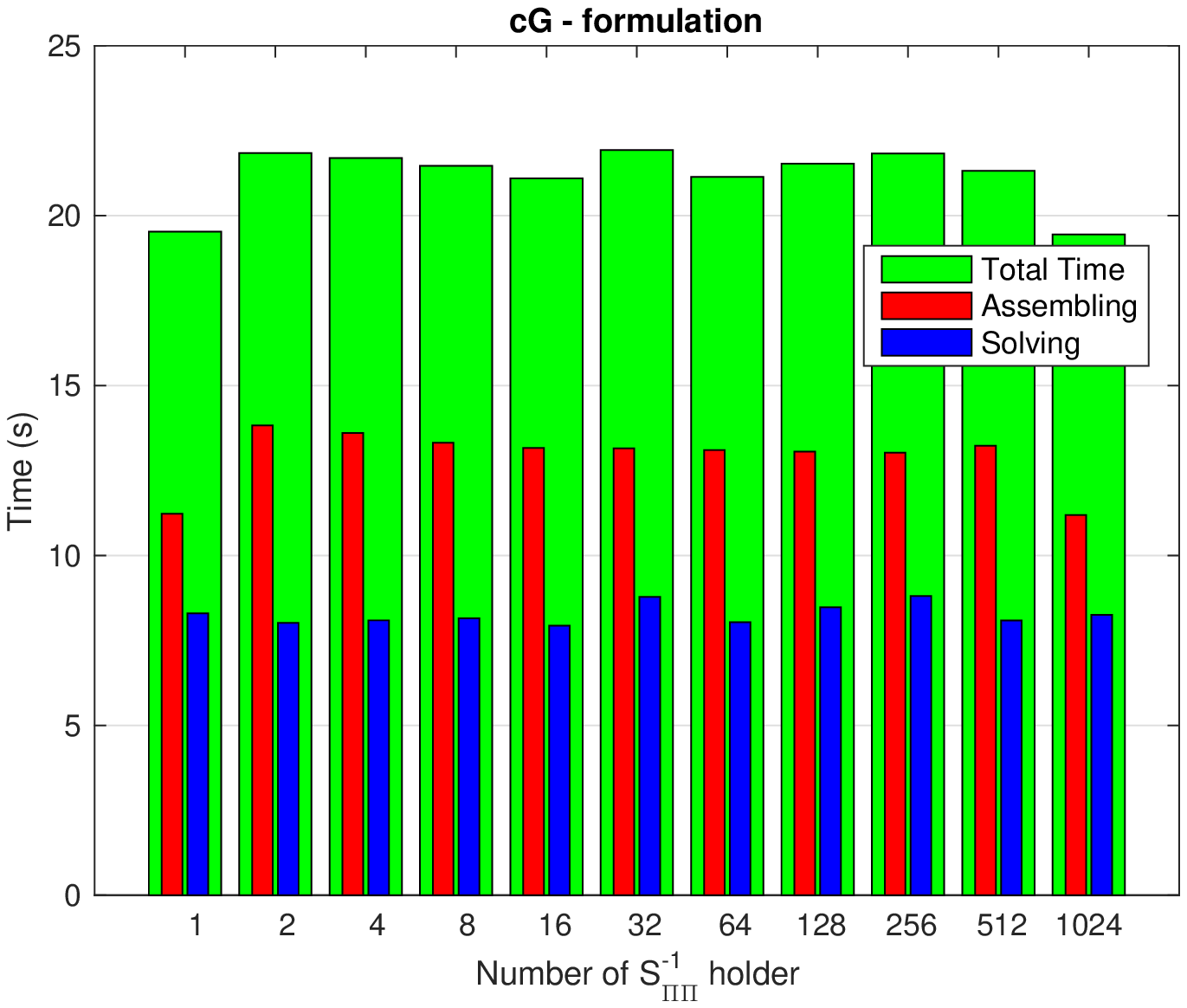}
  	\caption{$p=4$}
  \end{subfigure}
  \\
    \begin{subfigure}{0.328\textwidth}
        \includegraphics[width=0.99\textwidth]{\PathToPic 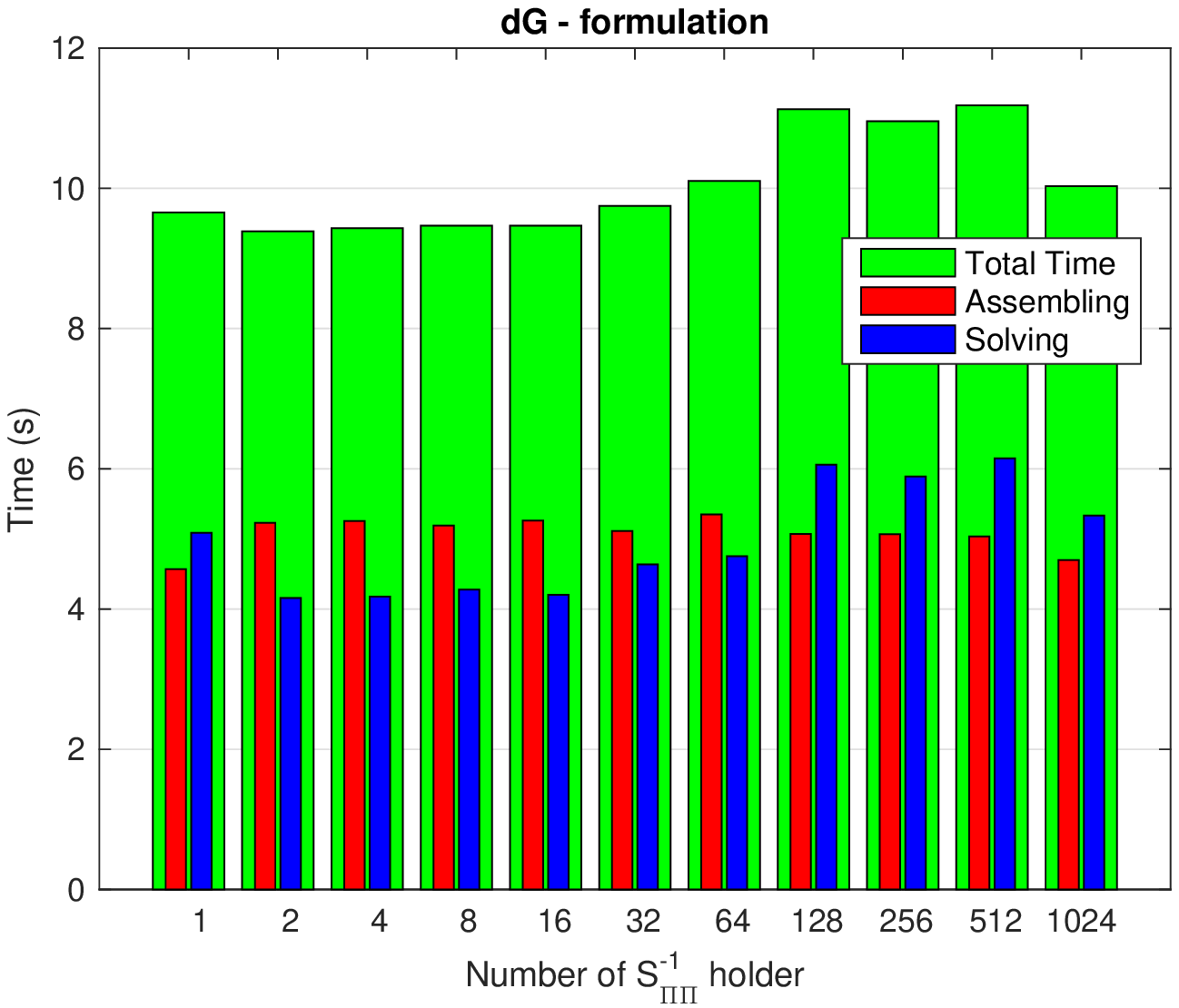}
  	\caption{$p=2$}
  \end{subfigure}
    \begin{subfigure}{0.328\textwidth}
        \includegraphics[width=0.99\textwidth]{\PathToPic 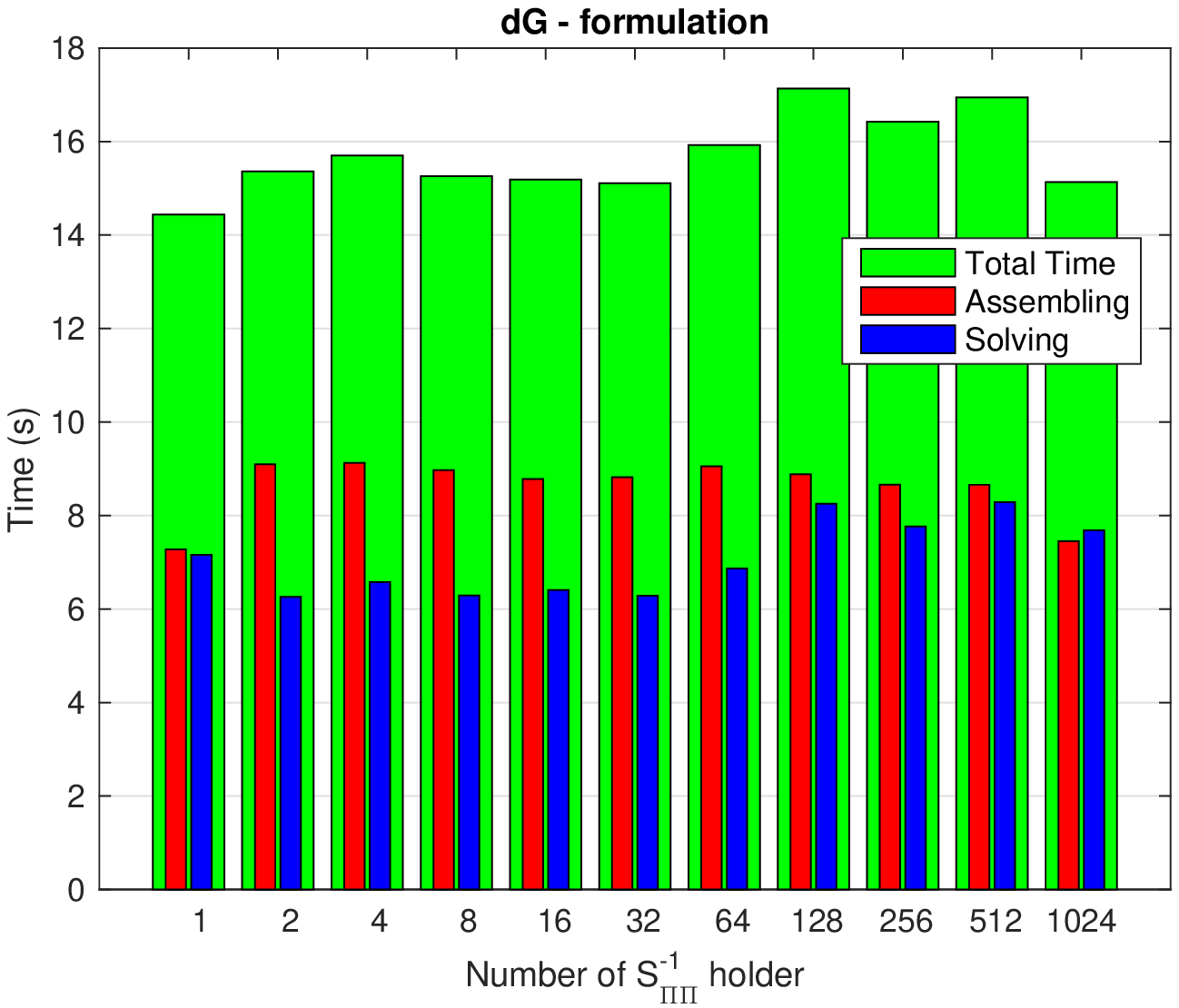}
  	\caption{$p=3$}
  \end{subfigure}
    \begin{subfigure}{0.328\textwidth}
        \includegraphics[width=0.99\textwidth]{\PathToPic 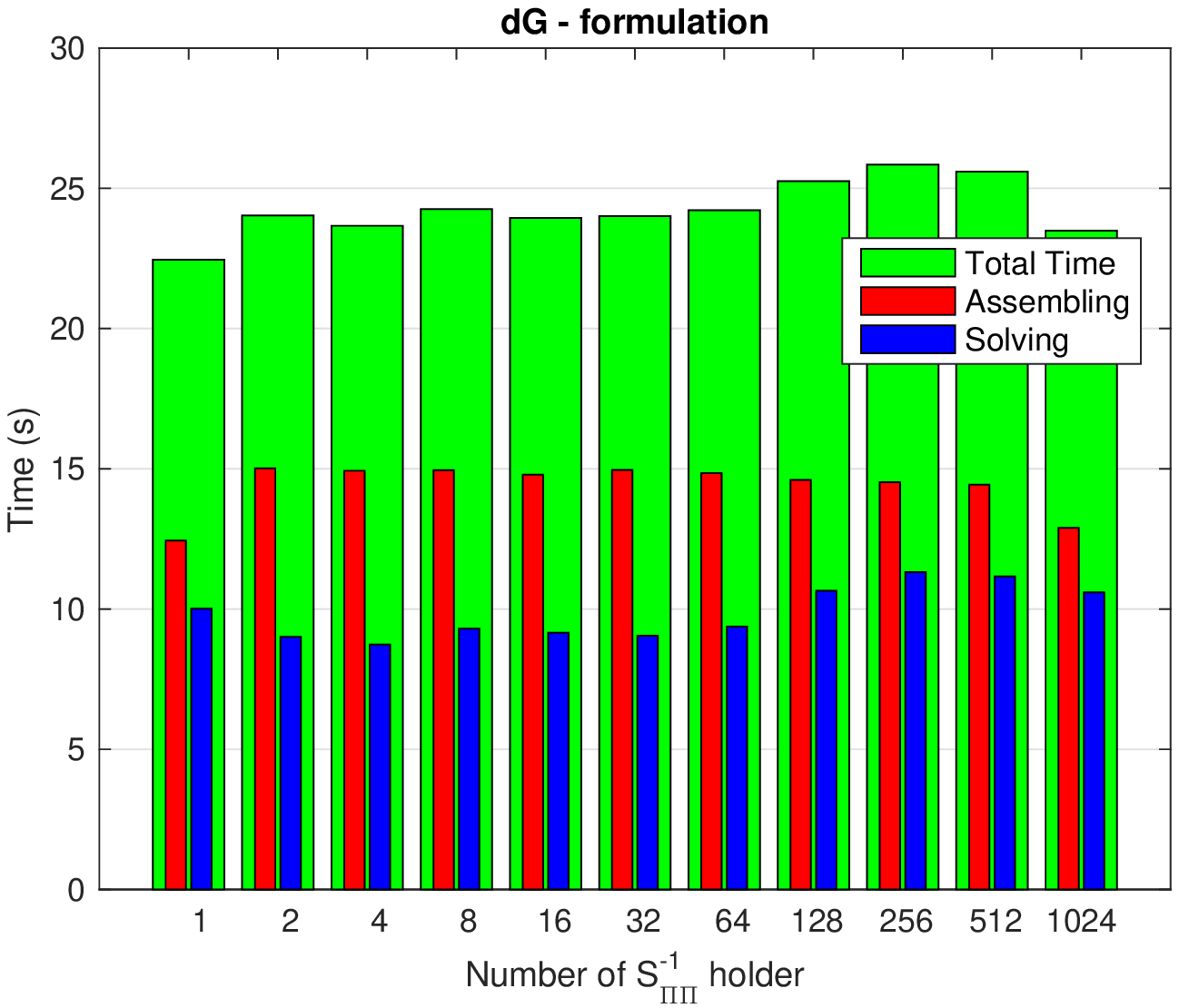}
  	\caption{$p=4$}
  \end{subfigure}
  \caption{Influence of the number of $S_{\Pi\Pi}^{-1}$ holders on the scaling. First row corresponds to cG-IETI-DP, second row to dG-IETI-DP. Each column has a fixed degree $p\in\{2,3,4\}$. Figures (a-c) summarizes the cG version and Figures (d-f) the dG version, respectively.}
  \label{fig:diffHolder}
\end{figure}

We only consider the two dimensional domain, where we perform $7$ initials refinements, but on a decomposition with $4096$ subdomains and end up with  around $70$~Mio. dofs. This gives a comparable setting as in Section~\ref{sec:weak} having the most refined domain. In order to better observe the influence of the number of $S_{\Pi\Pi}^{-1}$ holders, we increase the number of subdomains, leading to a larger coarse grid problem. We only investigate the case of using 1024 processors and the number of $S_{\Pi\Pi}$ holders given by $2^n,n\in\{0,1,\ldots,10\}$. Hence, we obtain the number of master processors ranging from 1 to 1024, such that each master has the same number of slaves. The results are summarized in Figure~\ref{fig:diffHolder} and Table~\ref{tab:diffHolder}.

\begin{table} [h]
      \begin{footnotesize}
 \begin{tabular}{|r|c|c|c||c|c|c||c|c|c|}\hline
\tiny{cG-IETI-DP}  & \multicolumn{3}{c||}{ $p=2$} & \multicolumn{3}{c||}{$p=3$}& \multicolumn{3}{c|}{$p=4$}\\ \hline
\makecell{$\#$ $S_{\Pi\Pi}^{-1}$ \\ Holder}& \makecell{Assemble \\ Time} & \makecell{Solving \\ Time}  & \makecell{Total \\ Time}   & \makecell{Assemble \\ Time} & \makecell{Solving \\ Time}  & \makecell{Total \\ Time}  & \makecell{Assemble \\ Time} & \makecell{Solving \\ Time}  & \makecell{Total \\ Time}      \\ \hline
  1    &   3.61 &   3.66 &  7.27 &   6.50  &   5.52 &12.02   &     11.23 &  8.30  &  19.53                      \\ \hline                                       
  2    &   4.49 &   3.58 &  8.07 &   7.97  &   5.57 &13.54   &     13.83 &  8.02  &  21.85                      \\ \hline                                       
  4    &   4.53 &   3.82 &  8.35 &   7.65  &   5.40 &13.05   &     13.60 &  8.09  &  21.69                      \\ \hline                                       
  8    &   4.46 &   3.63 &  8.09 &   7.72  &   5.76 &13.48   &     13.32 &  8.15  &  21.47                      \\ \hline                                       
  16   &   4.34 &   3.49 &  7.83 &   7.64  &   5.61 &13.25   &     13.16 &  7.93  &  21.09                      \\ \hline                                       
  32   &   4.33 &   3.73 &  8.06 &   7.74  &   5.39 &13.13   &     13.15 &  8.78  &  21.93                      \\ \hline                                       
  64   &   4.34 &   3.59 &  7.93 &   7.62  &   5.45 &13.07   &     13.10 &  8.04  &  21.14                      \\ \hline                                       
  128  &   4.49 &   4.06 &  8.55 &   7.60  &   6.05 &13.65   &     13.06 &  8.47  &  21.53                      \\ \hline                                       
  256  &   4.31 &   4.64 &  8.95 &   7.63  &   6.43 &14.06   &     13.02 &  8.81  &  21.83                      \\ \hline                                       
  512  &   4.34 &   3.61 &  7.95 &   7.55  &   5.71 &13.26   &     13.23 &  8.09  &  21.32                      \\ \hline                                       
  1024 &   3.73 &   3.80 &  7.53 &   6.56  &   5.77 &12.33   &     11.19 &  8.26  &  19.45                      \\ \hline\hline  
  \tiny{dG-IETI-DP}  & \multicolumn{3}{c||}{ $p=2$} & \multicolumn{3}{c||}{$p=3$}& \multicolumn{3}{c|}{$p=4$}\\ \hline
\makecell{$\#$ $S_{\Pi\Pi}^{-1}$ \\ Holder}  & \makecell{Assemble \\ Time} & \makecell{Solving \\ Time}  & \makecell{Total \\ Time}    & \makecell{Assemble \\ Time} & \makecell{Solving \\ Time}  & \makecell{Total \\ Time}  & \makecell{Assemble \\ Time} & \makecell{Solving \\ Time}  & \makecell{Total \\ Time}     \\ \hline
  1     & 4.57&  5.09&   9.66& 7.28  &  7.16 &  14.44& 12.44 & 10.01 &  22.45                      \\ \hline                                            
  2     & 5.23&  4.16&   9.39& 9.10  &  6.26 &  15.36& 15.02 &  9.01 &  24.03                      \\ \hline                                             
   4    & 5.25&  4.18&   9.43& 9.12  &  6.58 &  15.70& 14.93 &  8.73 &  23.66        \\ \hline                                                           
   8    & 5.19&  4.28&   9.47& 8.97  &  6.29 &  15.26& 14.95 &  9.30 &  24.25        \\ \hline                                                           
 16     & 5.26&  4.20&   9.46& 8.78  &  6.41 &  15.19& 14.79 &  9.16 &  23.95        \\ \hline                                                           
 32     & 5.11&  4.64&   9.75& 8.82  &  6.29 &  15.11& 14.96 &  9.05 &  24.01        \\ \hline                                                           
 64     & 5.35&  4.75&   10.1& 9.06  &  6.87 &  15.93& 14.85 &  9.37 &  24.22        \\ \hline                                                           
 128    & 5.07&  6.06&  11.13& 8.88  &  8.25 &  17.13& 14.61 & 10.65 &  25.26        \\ \hline                                                           
 256    & 5.07&  5.89&  10.96& 8.66  &  7.77 &  16.43& 14.52 & 11.32 &  25.84        \\ \hline                                                           
 512    & 5.03&  6.15&  11.18& 8.66  &  8.29 &  16.95& 14.43 & 11.16 &  25.59        \\ \hline                                                           
 1024   & 4.70&  5.33&  10.03& 7.45  &  7.68 &  15.13& 12.89 & 10.60 &  23.49        \\ \hline                                                           
 \end{tabular}
  \end{footnotesize}
    \caption{Influence of the number of processors having an LU-factorization of $S_{\Pi\Pi}$. Timings in seconds for 1024 Processors on a domain with around $70$~Mio. dofs and 2048 subdomains.  }
    \label{tab:diffHolder}
\end{table}

We observe that choosing several holders of the coarse grid problem in the cG version does not really have a significant effect. However, in the dG version, due to an increased number of primal variables, the use of several holders actually increases the performance of the solver by around $10\%$. Nevertheless, what is gained in the solving part does not pay off with the additional effort in the assembling phase. Considering the total computation time in Table~\ref{tab:diffHolder}, the best options is still either using only a single coarse grid problem on one processor or making a redundant factorization on each processor.

\section{Conclusion}
\label{sec:conclusion}
We have investigated the parallel scalability of the cG-IETI-DP and dG-IETI-DP method, respectively. Numerical tests showed a very good scalability in the strong and weak scaling for the assembling phase for both methods. We reached a speedup of approximately 900 when using 1024 cores. Although the speedup of the solver phase is not as good as the one for the assembler phase, we still reached a speedup of around 500 when using 1024 cores. One can even increase the parallel performance of the solver part by increasing the number of processors, which are holding the coarse grid problem. However, numerical examples have shown that this does not really pay off in the total time, due to an increased assembling time. To summarize, we saw that the proposed methods are well suited for large scale parallelization of assembling and solving IgA equations in two and three dimensions.

\section*{Acknowledgements}
This work was supported by the Austrian Science Fund (FWF) under the grant W1214, project DK4. This support is gratefully acknowledged. Moreover, the author wants to thank Prof. Ulrich Langer for the valuable comments and support during the preparation of the paper. He also gratefully acknowledges the support of Katharina Rafetseder of the Institute of Numerical Mathematics, Johannes Kepler University Linz (JKU) and Ioannis Toulopoulos and Angelos Mantzaflaris, Radon Institute of Computational and Applied Mathematics Linz (RICAM).

\section*{References}
\bibliographystyle{abbrv}

\bibliography{IETI_Parallelization.bib}

\begin{thebibliography}{10}

\bibitem{HL:BazilevsVeigaCottrellHughesSangalli:2006a}
Y.~{Bazilevs}, L.~{Beir\~ao da Veiga}, J.~A. {Cottrell}, T.~J.~R. {Hughes}, and
  G.~{Sangalli}.
\newblock {Isogeometric analysis: Approximation, stability and error estimates
  for $h$-refined meshes.}
\newblock {\em {Math. Models Methods Appl. Sci.}}, 16(7):1031--1090, 2006.

\bibitem{HL:BazilevsCaloCottrellEvans:2010a}
Y.~Bazilevs, V.~Calo, J.~Cottrell, J.~Evans, T.~Hughes, S.~Lipton, M.~Scott,
  and T.~Sederberg.
\newblock Isogeometric analysis using {T}-splines.
\newblock {\em Computer Methods in Applied Mechanics and Engineering},
  199(5–8):229 -- 263, 2010.
\newblock Computational Geometry and Analysis.

\bibitem{HL:BeiraodaVeigaBuffaSangalliVazquez:2014a}
L.~{Beir\~ao da Veiga}, A.~Buffa, G.~Sangalli, and R.~V\'{a}zquez.
\newblock Mathematical analysis of variational isogeometric methods.
\newblock {\em Acta Numerica}, 23:157--287, 2014.

\bibitem{HL:VeigaChoPavarinoScacchi:2012a}
L.~{Beir\~ao da Veiga}, D.~{Cho}, L.~F. {Pavarino}, and S.~{Scacchi}.
\newblock {Overlapping Schwarz methods for isogeometric analysis.}
\newblock {\em {SIAM J. Numer. Anal.}}, 50(3):1394--1416, 2012.

\bibitem{HL:VeigaChoPavarinoScacchi:2013a}
L.~{Beir\~ao Da Veiga}, D.~{Cho}, L.~F. {Pavarino}, and S.~{Scacchi}.
\newblock {BDDC preconditioners for isogeometric analysis.}
\newblock {\em {Math. Models Methods Appl. Sci.}}, 23(6):1099--1142, 2013.

\bibitem{HL:VeigaChoPavarinoScacchi:2013b}
L.~{Beir\~ao da Veiga}, D.~{Cho}, L.~F. {Pavarino}, and S.~{Scacchi}.
\newblock {Isogeometric Schwarz preconditioners for linear elasticity systems.}
\newblock {\em {Comput. Methods Appl. Mech. Eng.}}, 253:439--454, 2013.

\bibitem{HL:CotrellHughesBazilevs:2009a}
J.~A. Cotrell, T.~J.~R. Hughes, and Y.~Bazilevs.
\newblock {\em Isogeometric {A}nalysis, Toward {I}ntegration of {CAD} and
  {FEA}}.
\newblock John Wiley and Sons, 2009.

\bibitem{HL:PietroErn:2012a}
D.~A. {Di Pietro} and A.~{Ern}.
\newblock {\em {Mathematical aspects of discontinuous Galerkin methods.}}
\newblock Berlin: Springer, 2012.

\bibitem{HL:Dohrmann:2003a}
C.~R. {Dohrmann}.
\newblock {A preconditioner for substructuring based on constrained energy
  minimization.}
\newblock {\em {SIAM J. Sci. Comput.}}, 25(1):246--258, 2003.

\bibitem{HL:Dohrmann:2007a}
C.~R. Dohrmann.
\newblock An approximate {BDDC} preconditioner.
\newblock {\em Numerical Linear Algebra with Applications}, 14(2):149--168,
  2007.

\bibitem{HL:DouglasHaaseLanger:2003a}
C.~C. Douglas, G.~Haase, and U.~Langer.
\newblock {\em Tutorial on Elliptic PDE Solvers and Their Parallelization}.
\newblock Society for Industrial and Applied Mathematics, Philadelphia, PA,
  USA, 2003.

\bibitem{HL:DryjaGalvisSarkis:2007a}
M.~{Dryja}, J.~{Galvis}, and M.~{Sarkis}.
\newblock {BDDC methods for discontinuous Galerkin discretization of elliptic
  problems.}
\newblock {\em {J. Complexity}}, 23(4-6):715--739, 2007.

\bibitem{HL:DryjaGalvisSarkis:2013a}
M.~{Dryja}, J.~{Galvis}, and M.~{Sarkis}.
\newblock {A FETI-DP preconditioner for a composite finite element and
  discontinuous Galerkin method.}
\newblock {\em {SIAM J. Numer. Anal.}}, 51(1):400--422, 2013.

\bibitem{HL:DryjaSarkis:2014a}
M.~Dryja and M.~Sarkis.
\newblock 3-d feti-dp preconditioners for composite finite
  element-discontinuous galerkin methods.
\newblock In {\em Domain Decomposition Methods in Science and Engineering XXI},
  pages 127--140. Springer, 2014.

\bibitem{HL:GiannelliJuettlerSpeleers:2012a}
C.~Giannelli, B.~J\"uttler, and H.~Speleers.
\newblock {THB}-splines: the truncated basis for hierarchical splines.
\newblock {\em Comput. Aided Geom. Design}, 29, 2012.

\bibitem{HL:GiannelliJuettlerSpeleers:2014a}
C.~Giannelli, B.~J\"uttler, and H.~Speleers.
\newblock Strongly stable bases for adaptively refined multilevel spline
  spaces.
\newblock {\em Advances in Computational Mathematics}, 40:459--490, 2014.

\bibitem{HL:eigenweb}
G.~Guennebaud, B.~Jacob, et~al.
\newblock Eigen v3.
\newblock http://eigen.tuxfamily.org, 2010.

\bibitem{HL:Hofer:2016a}
C.~Hofer.
\newblock Analysis of discontinuous {G}alerkin dual-primal isogeometric tearing
  and interconnecting methods.
\newblock Technical Report No. 2016-03, {\
  }https://www.dk-compmath.jku.at/publications/dk-reports/2016-11-03, DK
  Computational Mathematics Linz Report Series, 2016.

\bibitem{HL:HoferLanger:2016b}
C.~Hofer and U.~Langer.
\newblock Dual-primal isogeometric tearing and interconnecting methods.
\newblock In P.~Neittanmakki, J.~Periaux, and O.~Pironneau, editors, {\em
  Contributions to PDE for Applications}, Springer-ECCOMAS series
  ``Computational Methods in Applied Sciences''. Springer, Berlin, Heidelberg,
  New York, 2016.
\newblock to appear.

\bibitem{HL:HoferLanger:2016a}
C.~Hofer and U.~Langer.
\newblock Dual-primal isogeometric tearing and interconnecting solvers for
  multipatch d{G}-{I}g{A} equations.
\newblock {\em Computer Methods in Applied Mechanics and Engineering}, 2016.
\newblock In Press, Accepted Manuscript,
  http://dx.doi.org/10.1016/j.cma.2016.03.031.

\bibitem{HL:HoferLangerToulopoulos:2016a}
C.~Hofer, U.~Langer, and I.~Toulopoulos.
\newblock {Discontinuous Galerkin Isogeometric Analysis of Elliptic Diffusion
  Problems on Segmentations with Gaps}.
\newblock {\em SIAM Journal on Scientific Computing}, 38(6):A3430--A3460, 2016.
\newblock available also at: http://arxiv.org/abs/1511.05715.

\bibitem{HL:HoferLangerToulopoulos:2016b}
C.~Hofer, U.~Langer, and I.~Toulopoulos.
\newblock Discontinuous {G}alerkin {I}sogeometric {A}nalysis on non-matching
  segmentation: error estimates and efficient solvers.
\newblock RICAM-Report~23, Johann Radon Institute for Computational and Applied
  Mathematics, Austrian Academy of Sciences, 2016.
\newblock available at
  https://www.ricam.oeaw.ac.at/publications/ricam-reports/Report No. 2016-23.

\bibitem{HL:HoferToulopoulos:2016a}
C.~Hofer and I.~Toulopoulos.
\newblock {Discontinuous Galerkin Isogeometric Analysis of elliptic problems on
  segmentations with non-matching interfaces}.
\newblock {\em Computers \& Mathematics with Applications}, 72(7):1811 -- 1827,
  2016.

\bibitem{Hoschek_Lasser_CAD_book_1993}
J.~Hoschek and D.~Lasser.
\newblock {\em Fundamentals of {C}omputet {A}ided {G}eometric {D}esign}.
\newblock A K Peters, Wellesley, Massachusetts, 1993.
\newblock Translated by L. Schumaker.

\bibitem{HL:HughesCottrellBazilevs:2005a}
T.~J.~R. Hughes, J.~A. Cottrell, and Y.~Bazilevs.
\newblock Isogeometric analysis: {CAD}, finite elements, {NURBS}, exact
  geometry and mesh refinement.
\newblock {\em Comput. Methods Appl. Mech. Engrg.}, 194:4135--4195, 2005.

\bibitem{HL:JuettlerKaplNguyenPanPauley:2014a}
B.~J\"uttler, M.~Kapl, D.-M. Nguyen, Q.~Pan, and M.~Pauley.
\newblock Isogeometric segmentation: The case of contractible solids without
  non-convex edges.
\newblock {\em Computer-Aided Design}, 57:74--90, 2014.

\bibitem{HL:KlawonnLanserRheinbach:2014a}
A.~{Klawonn}, M.~{Lanser}, and O.~{Rheinbach}.
\newblock {Nonlinear FETI-DP and BDDC methods.}
\newblock {\em {SIAM J. Sci. Comput.}}, 36(2):737--765, 2014.

\bibitem{HL:KlawonnLanserRheinbach:2015a}
A.~{Klawonn}, M.~{Lanser}, and O.~{Rheinbach}.
\newblock {Toward extremely scalable nonlinear domain decomposition methods for
  elliptic partial differential equations.}
\newblock {\em {SIAM J. Sci. Comput.}}, 37(6):c667--c696, 2015.

\bibitem{HL:KlawonnLanserRheinbach:2016a}
A.~{Klawonn}, M.~{Lanser}, and O.~{Rheinbach}.
\newblock {A nonlinear FETI-DP method with an inexact coarse problem.}
\newblock In {\em {Domain decomposition methods in science and engineering
  XXII. Proceedings of the 22nd international conference on domain
  decomposition methods, Lugano, Switzerland, September 16--20, 2013}}, pages
  41--52. Cham: Springer, 2016.

\bibitem{HL:KlawonnLanserRheinbachStengelWellein:2015a}
A.~Klawonn, M.~Lanser, O.~Rheinbach, H.~Stengel, and G.~Wellein.
\newblock {\em Hybrid MPI/OpenMP Parallelization in FETI-DP Methods}, pages
  67--84.
\newblock Springer International Publishing, Cham, 2015.

\bibitem{HL:KlawonnRheinbach:2006a}
A.~{Klawonn} and O.~{Rheinbach}.
\newblock A parallel implementation of {D}ual-{P}rimal {FETI} methods for
  three-dimensional linear elasticity using a transformation of basis.
\newblock {\em SIAM Journal on Scientific Computing}, 28(5):1886--1906, 2006.

\bibitem{HL:KlawonnRheinbach:2007a}
A.~Klawonn and O.~Rheinbach.
\newblock Inexact {FETI-DP} methods.
\newblock {\em International journal for numerical methods in engineering},
  69(2):284--307, 2007.

\bibitem{HL:KlawonnRheinbach:2010a}
A.~Klawonn and O.~Rheinbach.
\newblock Highly scalable parallel domain decomposition methods with an
  application to biomechanics.
\newblock {\em ZAMM - Journal of Applied Mathematics and Mechanics /
  Zeitschrift für Angewandte Mathematik und Mechanik}, 90(1):5--32, 2010.

\bibitem{HL:KlawonnRheinbachPavarino:2008a}
A.~{Klawonn}, O.~{Rheinbach}, and L.~F. {Pavarino}.
\newblock {Exact and inexact FETI-DP methods for spectral elements in two
  dimensions.}
\newblock In {\em {Domain decomposition methods in science and engineering
  XVII. Selected papers based on the presentations at the 17th international
  conference on domain decomposition methods, St. Wolfgang/Strobl, Austria,
  July 3--7, 2006.}}, pages 279--286. Berlin: Springer, 2008.

\bibitem{HL:KleissPechsteinJuettlerTomar:2012a}
S.~Kleiss, C.~Pechstein, B.~J{\"u}ttler, and S.~Tomar.
\newblock {IETI}--isogeometric tearing and interconnecting.
\newblock {\em Computer Methods in Applied Mechanics and Engineering},
  247:201--215, 2012.

\bibitem{HL:PARDISO500}
A.~Kuzmin, M.~Luisier, and O.~Schenk.
\newblock Fast methods for computing selected elements of the greens function
  in massively parallel nanoelectronic device simulations.
\newblock In F.~Wolf, B.~Mohr, and D.~Mey, editors, {\em Euro-Par 2013 Parallel
  Processing}, volume 8097 of {\em Lecture Notes in Computer Science}, pages
  533--544. Springer Berlin Heidelberg, 2013.

\bibitem{HL:LangerMantzaflarisMooreToulopoulos:2015b}
U.~Langer, A.~Mantzaflaris, S.~E. Moore, and I.~Toulopoulos.
\newblock Multipatch discontinuous {G}alerkin isogeometric analysis.
\newblock In B.~J\"uttler and B.~Simeon, editors, {\em Isogeometric Analysis
  and Applications IGAA 2014}, volume 107 of {\em Lecture Notes in Computer
  Science}, pages 1--32, Heidelberg, 2015. Springer.
\newblock also available at http://arxiv.org/abs/1411.2478.

\bibitem{HL:LangerToulopoulos:2015a}
U.~Langer and I.~Toulopoulos.
\newblock {Analysis of multipatch discontinuous Galerkin IgA approximations to
  elliptic boundary value problems}.
\newblock {\em Computing and Visualization in Science}, 17(5):217--233, 2015.

\bibitem{HL:LiWidlund:2007a}
J.~{Li} and O.~B. {Widlund}.
\newblock {On the use of inexact subdomain solvers for BDDC algorithms.}
\newblock {\em {Comput. Methods Appl. Mech. Eng.}}, 196(8):1415--1428, 2007.

\bibitem{HL:MandelDohrmannTezaur:2005a}
J.~{Mandel}, C.~R. {Dohrmann}, and R.~{Tezaur}.
\newblock {An algebraic theory for primal and dual substructuring methods by
  constraints.}
\newblock {\em {Appl. Numer. Math.}}, 54(2):167--193, 2005.

\bibitem{gismoweb}
A.~Mantzaflaris, C.~Hofer, et~al.
\newblock {G+Smo (Geometry plus Simulation modules) v0.8.1}.
\newblock http://gs.jku.at/gismo, 2015.

\bibitem{HL:PauleyNguyenMayerSpehWeegerJuettler:2015a}
M.~Pauley, D.-M. Nguyen, D.~Mayer, J.~Speh, O.~Weeger, and B.~J\"uttler.
\newblock The isogeometric segmentation pipeline.
\newblock In B.~J\"uttler and B.~Simeon, editors, {\em Isogeometric Analysis
  and Applications IGAA 2014}, Lecture Notes in Computer Science, Heidelberg,
  2015. Springer.
\newblock to appear, also available as Technical Report no. 31 at
  http://www.gs.jku.at.

\bibitem{HL:Pechstein:2013a}
C.~{Pechstein}.
\newblock {\em {Finite and boundary element tearing and interconnecting solvers
  for multiscale problems.}}
\newblock Berlin: Springer, 2013.

\bibitem{HL:Rheinbach:2009}
O.~{Rheinbach}.
\newblock {Parallel iterative substructuring in structural mechanics.}
\newblock {\em {Arch. Comput. Methods Eng.}}, 16(4):425--463, 2009.

\bibitem{HL:Riviere:2008a}
B.~{Rivi\`ere}.
\newblock {\em {Discontinuous Galerkin methods for solving elliptic and
  parabolic equations. Theory and implementation.}}
\newblock Philadelphia, PA: Society for Industrial and Applied Mathematics
  (SIAM), 2008.

\bibitem{HL:ToselliWidlund:2005a}
A.~{Toselli} and O.~B. {Widlund}.
\newblock {\em {Domain decomposition methods -- algorithms and theory.}}
\newblock Berlin: Springer, 2005.

\bibitem{HL:Tu:2007b}
X.~Tu.
\newblock Three-level {BDDC} in three dimensions.
\newblock {\em SIAM Journal on Scientific Computing}, 29(4):1759--1780, 2007.

\bibitem{HL:Tu:2007a}
X.~Tu.
\newblock Three-level {BDDC} in two dimensions.
\newblock {\em International journal for numerical methods in engineering},
  69(1):33--59, 2007.

\bibitem{HL:Zampini:2014}
S.~Zampini.
\newblock {\em Inexact {BDDC} Methods for the Cardiac Bidomain Model}, pages
  247--255.
\newblock Springer International Publishing, Cham, 2014.

\end{thebibliography}

\end{document}